\documentclass{article}
\usepackage{latexsym,amsfonts,amsmath,amsthm,makeidx}
\usepackage{makeidx}

\makeindex
\newtheorem{definition}{Definition}[section]
\newtheorem{theorem}{Theorem}[section]
\newtheorem{lemma}{Lemma}[section]
\newtheorem{corollary}{Corollary}[section]
\newtheorem{proposition}{Proposition}[section]
\newtheorem{remark}{Remark}[section]
\newtheorem{example}{Example}[section]

\newcommand{\s}{\section}

\newcommand{\R}{\mathbb R}

\newcommand{\lab}{\label}
\newcommand{\bt}{\begin{theorem}}
\newcommand{\et}{\end{theorem}}
\newcommand{\bl}{\begin{lemma}}
\newcommand{\el}{\end{lemma}}
\newcommand{\bd}{\begin{definition}}
\newcommand{\ed}{\end{definition}}
\newcommand{\bc}{\begin{corollary}}
\newcommand{\ec}{\end{corollary}}
\newcommand{\bp}{\begin{proof}}
\newcommand{\ep}{\end{proof}}
\newcommand{\bx}{\begin{example}}
\newcommand{\ex}{\end{example}}
\newcommand{\bi}{\begin{exercise}}
\newcommand{\ei}{\end{exercise}}
\newcommand{\bo}{\begin{proposition}}
\newcommand{\eo}{\end{proposition}}
\newcommand{\br}{\begin{remark}}
\newcommand{\er}{\end{remark}}
\newcommand{\be}{\begin{equation}}
\newcommand{\ee}{\end{equation}}
\newcommand{\ba}{\begin{align}}
\newcommand{\ea}{\end{align}}
\newcommand{\bn}{\begin{enumerate}}
\newcommand{\en}{\end{enumerate}}
\newcommand{\bg}{\begin{align*}}
\newcommand{\bcs}{\begin{cases}}
\newcommand{\ecs}{\end{cases}}

\newcommand{\bean}{\begin{eqnarray*}}
\newcommand{\eean}{\end{eqnarray*}}

\numberwithin{equation}{section}

\title{\bf Multiplicity and concentration behavior of solutions to  the critical Kirchhoff type problem\thanks{Supported
by NSFC(11401583) and the Fundamental Research Funds for the central Universities
(16CX02051A).  } }
\author{{\bf \small Jian Zhang}\\
 \small\it College of Science, China University of Petroleum,\\
\small\it Qingdao 266580, Shandong, P. R. China\\
\small\it zjian@upc.edu.cn\\
\\
{\bf \small Wenming  Zou}\\
 \small\it Department of Mathematical Sciences, Tsinghua University,\\
 \small\it Beijing 100084, P. R. China\\
\small\it wzou@math.tsinghua.edu.cn
}
\date{}
\begin{document}
\maketitle
 \pagenumbering{arabic} \thispagestyle{empty}
\setcounter{page}{1} \vspace{0.2cm} {\footnotesize \noindent
 \line(1,0){345}\\[10pt]
 {\bf Abstract} \
In this paper, we study the multiplicity and concentration of the positive solutions to  the following critical Kirchhoff type problem:
\begin{equation*}
 -\left(\varepsilon^2 a+\varepsilon b\int_{\R^3}|\nabla u|^2\mathrm{d} x\right)\Delta u + V(x) u  = f(u)+u^5\ \
{\rm in } \ \  \R^3,
\end{equation*}
where $\varepsilon$ is a small positive parameter, $a$, $b$ are positive constants, $V \in C(\mathbb{R}^3)$ is a positive potential,
$f \in C^1(\R^+, \R)$ is a subcritical nonlinear term, $u^5$ is a pure critical nonlinearity. When $\varepsilon>0$ small, we establish the relationship between the number of positive solutions and the profile of the potential $V$. The exponential decay at infinity of the solution is also obtained. In particular,    we show that  each solution concentrates around a local strict minima  of $V$ as $\varepsilon \rightarrow 0$.
\vspace{0.2cm} 
}\\[2mm]
 \line(1,0){345}

\date{}

\numberwithin{equation}{section}

\maketitle

\newpage

\s{Introduction}
\renewcommand{\theequation}{1.\arabic{equation}}

This paper is concerned with the existence and concentration of multiple positive solutions to  the following critical Kirchhoff type equation:
\begin{equation} \label{1.1}
 -\left(\varepsilon^2 a+\varepsilon b\int_{\R^3}|\nabla u|^2\mathrm{d} x\right)\Delta u + V(x) u  = f(u)+u^5\ \
{\rm in } \ \  \R^3,
\end{equation}
where $\varepsilon$ is a small positive parameter, $a$, $b$ are positive constants, $V$ and $f$ are continuous functions satisfying some additional assumptions.

\vskip0.11in
The Kirchhoff equation  occurs in various branches of mathematical
physics. For example, it can be used to model suspension bridges (see \cite{A-B-G}). In  particular,  the following    problem
\begin{equation}\label{1.2}
 -\left(a+b\int_{\Omega}|\nabla u|^2\mathrm{d} x\right)\Delta u  = f(x,u)\ \
{\rm in } \ \ \  \ \R^3,\\[2mm]
\end{equation}
is closely  related to the  wave equation  counterpart
\begin{equation*}
\rho \frac{\partial^2 u}{\partial t^2}-\left(\frac{P_0}{h}+\frac{E}{2L}\int_0^L\left|\frac{\partial u}{\partial x}\right|^2\mathrm{d}x\right)\frac{\partial^2 u}{\partial x^2}=0
\end{equation*}
and

\begin{equation}\label{1.2-zwm}
u_{tt} -\left(a+b\int_{\Omega}|\nabla u|^2\mathrm{d} x\right)\Delta u  = f(x,u)
\end{equation}
proposed by Kirchhoff \cite{K}. Because of the presence of the nonlocal term $\left(\int_{\Omega}|\nabla u|^2\mathrm{d} x\right)\Delta u$, Eq. (\ref{1.2}) is not a pointwise identity, which causes additional mathematical difficulties.
For example, it is more difficult to check the geometric structure of the functional associated with the equation and the boundedness, convergence of the Palais-Smale sequence if we seek solutions using variational methods.
  Moreover, we cannot derive
\begin{align}\label{1-1}
\int_{\R^N}|\nabla u_n|^2 \mathrm{d}x \int_{\R^N}\nabla u_n \nabla v \mathrm{d}x \rightarrow \int_{\R^N}|\nabla u|^2 \mathrm{d}x \int_{\R^N} \nabla u \nabla v\mathrm{d}x, \ \  v \in H^1(\R^N)
\end{align}
and
\begin{align}\label{1-2}
\left(\int_{\R^N}|\nabla u_n|^2 \mathrm{d}x\right)^2 - \left(\int_{\R^N}|\nabla u|^2 \mathrm{d}x\right)^2= \left(\int_{\R^N}|\nabla u_n- \nabla u|^2 \mathrm{d}x\right)^2+o_n(1)
\end{align}
from $u_n \rightharpoonup u$ weakly in $H^1(\R^N)$, which is crucial when we consider the convergence of the Palais-Smale sequence.



 \vskip0.12in

When $b=0$, problem (\ref{1.1}) reduces to the singularly perturbed problem
\begin{equation}\label{1.4}
-\varepsilon ^2 \Delta u + V(x) u = g(u) \ \ \mathrm{in} \ \
\mathbb{R}^N.
\end{equation}
Many authors are concerned with problem (\ref{1.4}) for $\varepsilon>0$ small since solutions of (\ref{1.4}) are known as semiclassical states,
which can be used to describe the transition from quantum to classical mechanics. Let us recall some results. In $\cite{FW,O1,O2}$, the Lyapunov-Schmidt reduction is used to construct single and multiple spike solutions. However, the Lyapunov-Schmidt reduction is
based on the uniqueness or non-degeneracy of ground state solutions of the corresponding limiting equation. To overcome
this difficulty, Rabinowitz $\cite{Rabinowitz2}$ firstly used the variational approach to obtain the existence of solutions of (\ref{1.4}) for $\varepsilon>0$ small under the assumption
\begin{align}\label{1.5}
\liminf_{|x| \rightarrow \infty}V(x)> \inf_{x \in \R^N}V(x)>0.
\end{align}
In \cite{X-W}, the concentration of solutions was also proved. By introducing a penalization approach, a localized version of the result in $\cite{Rabinowitz2,X-W}$ was proved by del Pino-Felmer $\cite{DF}$.
Subsequently, Jeanjean-Tanaka $\cite{JT}$ extended $\cite{DF}$'s work to
a more general form. For other results on
the singularly perturbed problem, see \cite{ABC,ABC2,BJ,DF2,DF3,DF4,Gui} and the references therein.

\vskip0.2in

When $b \ne 0$, He-Zou \cite{HZ} firstly  studied the Kirchhoff type problem
\begin{equation*}
 -\left(\varepsilon^2 a+\varepsilon b\int_{\R^3}|\nabla u|^2\mathrm{d} x\right)\Delta u + V(x) u  = f(u)\ \
{\rm in } \ \  \R^3,
\end{equation*}
where $V$ is a positive continuous function and $f$ is a subcritical term. By using the Ljusternik-Schnirelmann theory and relating the number of solutions with the topology of the set where $V$ attains its minimum, they proved  the multiplicity and concentration behavior of positive solutions under the assumption \eqref{1.5}. For the critical case, Wang {\it et al.} \cite{WTXZ} considered the problem
\begin{equation}\label{zwm=1}
 -\left(\varepsilon^2 a+\varepsilon b\int_{\R^3}|\nabla u|^2\mathrm{d} x\right)\Delta u + V(x) u  = \lambda f(u)+ |u|^4u\ \
{\rm in } \ \  \R^3,
\end{equation}
where $V(x)$ admits at least one minimum. For $\varepsilon>0$ small enough and $\lambda>0$ sufficiently large, they obtained the existence, concentration and some further properties of the positive ground state solution to Eq. \eqref{zwm=1}. They also extended the results of \cite{HZ} to the critical case.
We emphasize   that, in \cite{WTXZ}, the assumption (\ref{1.5}) is required and plays an indispensable role for the arguments of  \cite{WTXZ}. Moreover, in order to deal with the critical term, the parameter $\lambda$ is required to  be large.
When $V(x)$ is a locally H\"{o}lder continuous function satisfying $\inf_{x \in \R^3}V(x)>0$ and $\inf_{\wedge}V < \min_{\partial \wedge} V$ for some open bounded set $\wedge\in \R^3$,  He {\it et al.} \cite{HLP} got a solution of  Eq. \eqref{zwm=1} concentrating around a local minimum point of $V$ in $\wedge$.
With the Ljusternik-Schnirelmann theory, they also got multiple solutions by employing the topology of the set where $V$ attains its local minimum.
Later, the main results in \cite{HLP} was extended by \cite{H-L} for the case $f(u)=|u|^{p-2}u$ $(2< p \le 4)$ in  Eq. \eqref{zwm=1}. For other related results, the readers may see \cite{G-F,FIS,LG2} and the references therein.

\vskip0.12in

Motivated by the above results, in this paper, we will study the  multiplicity and concentration behavior of positive solutions to  the critical problem (\ref{1.1}). Our results are different from the results mentioned above.  Before stating  the main result, we introduce   the following hypotheses:

\begin{itemize}
\item [$(V_{1})$] $V(x) \in C(\mathbb{R}^3)$ and $V_0=:\inf_{x \in \R^3}V(x)>0$;
\item [$(V_{2})$] there exist $x^1$, $x^2$, \ldots $x^k$ in $\R^3$ such that each point of $x^1$, $x^2$, \ldots $x^k$ is a strict global minima of $V$
in $\R^3$;
\item [$(f_1)$] $f \in C^1(\R^+, \R)$ and $\lim_{u \rightarrow 0+}\frac{f(u)}{u^3}=\lim_{u \rightarrow +\infty}\frac{f(u)}{u^5}=0$;
\item [$(f_2)$] the function $\frac{f(u)}{u^3}$ is increasing for $u > 0$ and $\lim_{u \rightarrow +\infty}\frac{f(u)}{u^3} = +\infty$.
\end{itemize}

\vskip0.123in

\bt\lab{Theorem 1.1} Assume
$(V_1)$-$(V_2)$ and $(f_1)$-$(f_2)$ hold. Then there exists $\varepsilon_0>0$ such that (\ref{1.1}) admits at least $k$ different positive solutions $v_\varepsilon^i$, $i=1$, $2$, \ldots $k$ for $\varepsilon \in (0,\varepsilon_0)$. Moreover, each $v_\varepsilon^i$ possesses a maximum point $z_\varepsilon^i \in \R^3$ such that $V(z_\varepsilon^i) \rightarrow V(x^i)=V_0$ as $\varepsilon \rightarrow 0$. Besides, there exist $C_0^i$, $c_0^i>0$ satisfying
\begin{align*}
v_\varepsilon^i(x) \le C_0^i \exp\left(-c_0^i \frac{|x-z_\varepsilon^i|}{\varepsilon}\right)
\end{align*}
for $\varepsilon \in (0,\varepsilon_0)$ and $x \in \R^3$.
\et

\br\label{zwm=2}  The novelties of Theorem \ref{Theorem 1.1} are twofold:   The Rabinowitz type assumption  (\ref{1.5}) is removed; the  multiplicity and asymptotic behavior of the solutions is obtained. 
\er

\vskip0.12in
\br Some ideas of the current paper is inspired by \cite{CN}, where  Cao-Noussair studied the subcritical problem
\begin{equation}\label{1.7}
 -\Delta u + \mu u  = Q(x) |u|^{p-2}u\ \
{\rm in } \ \  \R^N,
\end{equation}
where $2<p<2^*$. They obtained the existence of both positive and nodal solutions  to  (\ref{1.7}) which is  affected by the shape of the graph of $Q(x)$. Subsequently, the idea of \cite{CN} is applied in \cite{LG2} to deal with the singularly perturbed critical Schr\"{o}dinger-Poisson equation.
We note that in \cite{LG2}, the strict inequality (\ref{1.5}) is used to estimate the energy level,
which is crucial for the proof of the relative compactness of the Palais-Smale sequence.
Unfortunately, such a method does not work for problem (\ref{1.1}) since (\ref{1.5}) does not hold.
Compared with  the works in \cite{CN,LG2}, another major difficulty in dealing with (\ref{1.1}) lies in the presence of the nonlocal term $\left(\int_{\R^3}|\nabla u|^2\mathrm{d} x\right)\Delta u$. Since (\ref{1-1})-(\ref{1-2}) do not hold in general, it is difficult to prove the compactness of the Palais-Smale sequence. Thus, we need a deeper understanding of the obstructions to the compactness.
By developing some techniques, we can estimate the Palais-Smale sequence carefully and solve the problem.
 \er

\vskip0.2in
Before closing this section, we say a few words on the equation (\ref{1.1}) with $\varepsilon=1$. It becomes
\begin{equation} \label{1.3}
-\left(a+b\int_{\R^N}|\nabla u|^2\mathrm{d} x\right)\Delta u + V(x)u  = f(x,u)\ \
{\rm in } \ \ \R^N
\end{equation}
and  receives much attention in recent years. When $f(x,u)=f(u)$ and $V(x)$ is radially symmetric, the  working  space $H_r^1(\R^N)$ is compactly embedded into  $L^p(\R^N) \ (2<p<2^*)$, the authors  of  \cite{LLS1} recovered the   compactness and proved the existence of solutions. When Eq. (\ref{1.3}) is non-radial,  the existence of solutions was
obtained in \cite{SW}. In particular, in \cite{AF}, Alves-Figueiredo considered the periodic Kirchhoff equation with critical growth and proved the existence of positive solutions. The authors of \cite{LY} got the existence of a positive ground state solution to Eq. (\ref{1.3}) by using a monotonicity trick and a new version of global compactness lemma. Recently, the result of  \cite{LY} studying Eq. (\ref{1.3}) was extended to the critical case by \cite{LG1}.


\vskip0.12in

The outline of this paper is as follows: in Section $2$, we
establish some  key lemmas; in Section $3$, we prove the existence of multiple solutions of (\ref{1.1});
the  Section $4$ is devoted to prove the concentration of solutions of (\ref{1.1}).

\vskip0.15in

\noindent{\bf Notations:}
\begin{itemize}
\item [$\bullet$] $\|u\|_s:=\big(\int_{\mathbb{R}^3}|u|^s\mathrm{d}x\big)^{\frac{1}{s}}$,
$1 \le s \le \infty$;
\item [$\bullet$] $H=H^1(\mathbb{R}^3)$ denotes the Hilbert space equipped with the
norm $\|u\|_{H}^2=\int_{\mathbb{R}^3}(|\nabla u|^2 + |u|^2)\mathrm{d}
x$, $D^{1,2}(\mathbb{R}^3)=\left\{u \in L^6(\mathbb{R}^3): \nabla u \in L^2(\mathbb{R}^3)\right\}$ denotes the Sobolev space equipped with the norm $\|u\|^2_{D^{1,2}}=\int_{\mathbb{R}^3}|\nabla u|^2\mathrm{d} x$;
\item [$\bullet$] $S:=\inf_{u\in D^{1,2} (\mathbb{R}^3) \setminus \{0\}}
\frac{\int_{\mathbb{R}^3} |\nabla u|^2 \mathrm{d}
x}{\big(\int_{\mathbb{R}^3}|u|^6\mathrm{d}
x\big)^{\frac{1}{3}}}$ denotes the best Sobolev constant;
\item [$\bullet$] $C$ denotes a positive constant (possibly different).
\end{itemize}

\s{Preliminary Lemmas}

\renewcommand{\theequation}{2.\arabic{equation}}
We assume $f(u)=0$ for $u \le 0$. Note that $\liminf_{|x| \rightarrow \infty}V(x) \ge \inf_{x \in \R^3}V(x)$.
When $\liminf_{|x| \rightarrow \infty}V(x)=+\infty$, Theorem \ref{Theorem 1.1} can be proved easier because the embedding
\begin{align*}
\left\{u \in H^1(\R^3): \int_{\R^3}\left(\varepsilon^2 |\nabla u|^2 + V(\varepsilon x) |u|^2\right)\mathrm{d}x<+\infty\right\} \hookrightarrow L^p(\R^3), \ \ 2<p<6
\end{align*}
is compact. Thus, we only consider the case $\liminf_{|x| \rightarrow \infty}V(x)<+\infty$.
Denote $V_\infty=:\liminf_{|x| \rightarrow \infty}V(x)$. Then $V_\infty \ge V_0$. By $(V_2)$, we have $V(x^i)=V_0$ for $i=1$, $2$, \ldots $k$. By $(f_1)$-$(f_2)$, we derive
\begin{align*}
\frac{1}{4}f(u)u - F(u) \ge 0, \ \ f'(u)u^2-3f(u)u \ge 0, \ \ \forall\ u \in \R,
\end{align*}
where $f(u)u \ge 0$, $F(u)=\int_0^u f(s) \mathrm{d}s \ge 0$. Moreover, the function
$\frac{1}{4}f(u)u-F(u)$ is increasing for $u \in \R$.

Make the change of variable $\varepsilon z=x$, we can rewrite (\ref{1.1}) as
\begin{equation} \label{2.1}
 -\left(a+ b\int_{\R^3}|\nabla u|^2\mathrm{d} x\right)\Delta u + V(\varepsilon x) u  = f(u)+u^5\ \
{\rm in } \ \  \R^3.
\end{equation}
For any fixed $\varepsilon>0$, let $H_\varepsilon=\left\{u \in H: \int_{\R^3}V(\varepsilon x) |u|^2\mathrm{d}x<\infty\right\}$ be the Hilbert space with the inner product
$(u,v)_\varepsilon=\left(\int_{\R^3}a \nabla u \nabla v+ V(\varepsilon x) u v\mathrm{d}x\right)$.
Then the norm of $H_\varepsilon$ is $\|u\|_\varepsilon=\left(\int_{\R^3}|\nabla u|^2 +V(\varepsilon x) |u|^2\mathrm{d}x\right)^{\frac{1}{2}}$.
The functional associated with (\ref{2.1}) is
\begin{align}\label{2.2}
I_\varepsilon(u) = &\frac{1}{2} \|u\|_\varepsilon^2+\frac{b}{4}\left(\int_{\R^3}|\nabla u|^2 \mathrm{d} x\right)^2- \int_{\R^3}F(u) \mathrm{d}x-\frac{1}{6}\int_{\R^3}|u|^6\mathrm{d}x,
\end{align}
where $u \in H_\varepsilon$. Moreover, for any $v \in H_\varepsilon$,
\begin{align*}
(I_\varepsilon'(u),v)=&\int_{\R^3}\left(a \nabla u \nabla v+V(\varepsilon x) u v\right)\mathrm{d}x + b \int_{\R^3}|\nabla u|^2 \mathrm{d} x \int_{\R^3}\nabla u \nabla v\mathrm{d} x\\
&-\int_{\R^3} f(u) v \mathrm{d} x-\int_{\R^3}u^5 v\mathrm{d}x.
\end{align*}
Clearly, the functional $I_\varepsilon : H_\varepsilon \mapsto \mathbb{R}$
is of class $C^{1}$ and critical points of $I_\varepsilon$ are weak solutions of (\ref{2.1}).
Let
$m_\varepsilon=\inf\{I_\varepsilon(u): u \in M_\varepsilon\}$ with $M_\varepsilon=\{u \in H_\varepsilon \setminus \{0\}:
(I_\varepsilon'(u),u)=0\}$.

Let $d>0$.
Then for $u \in H$, we can define
$\|u\|_d=\left(\int_{\mathbb{R}^3}(a|\nabla u|^2+d |u|^2) \mathrm{d} x\right)^{\frac{1}{2}}$.
Obviously, $\|.\|_d$ is an equivalent norm on $H$.
Consider the following equation:
\begin{equation} \label{2.3}
 -\left(a+ b\int_{\R^3}|\nabla u|^2\mathrm{d} x\right)\Delta u + d u  = f(u)+u^5\ \
{\rm in } \ \  \R^3.
\end{equation}
The functional associated with (\ref{2.3}) is
\begin{align*}
I_d(u) = &\frac{1}{2} \|u\|_d^2+\frac{b}{4}\left(\int_{\R^3}|\nabla u|^2 \mathrm{d} x\right)^2- \int_{\R^3}F(u) \mathrm{d}x-\frac{1}{6}\int_{\R^3}|u|^6\mathrm{d}x.
\end{align*}
Let
$m_d=\inf\{I_d(u): u \in M_d\}$ with $M_d=\{u \in H \setminus \{0\}:
(I_{d}'(u),u)=0\}$.

Recall that $S$ is attained by the functions
$\frac{\varepsilon^{\frac{1}{4}}}{\big(\varepsilon +
|x|^{2}\big)^{\frac{1}{2}}}$, where $\varepsilon>0$.
Let $u_{\varepsilon}(x) = \frac{\psi(x) \varepsilon^{\frac{1}{4}}}{\big(\varepsilon +
|x|^{2}\big)^{\frac{1}{2}}}$, where $\psi \in C_{0}^{\infty}
(B_{2r}(0))$ such that $\psi(x) = 1$ on
$B_{r}(0)$ and $0 \leq \psi(x) \leq 1$. From $\cite{Willem}$, we have the following results.

\bl\lab{Lemma 2.5} For $\varepsilon>0$ small, there holds
\begin{align}\label{2.6}
&\int_{\mathbb{R}^3} |\nabla u_{\varepsilon}|^{2}
\mathrm{d} x =  K_{1} + O
(\varepsilon^{\frac{1}{2}}),\ \ \ \int_{\mathbb{R}^3} |u_{\varepsilon}|^6 \mathrm{d}
x=K_{2} + O (\varepsilon^{\frac{3}{2}}),\notag\\
&\int_{\mathbb{R}^3} |u_{\varepsilon}|^2 \mathrm{d} x =O\left(\varepsilon^{\frac{1}{2}}\right),
\end{align}
where $S =\frac{K_{1}}{K_{2}^\frac{1}{3}}$.

\el
\vskip0.123in

By Lemma $2.1$, we can get the following result.

\vskip0.123in

\bl\lab{Lemma 2.6} Let $d>0$. Then
\begin{align*}
m_d<\hat{c}:=\frac{a}{3}\left\{\frac{bS^3+\sqrt{(bS^3)^2+4aS^3}}{2}\right\}+\frac{b}{12}\left\{\frac{bS^3+\sqrt{(bS^3)^2+4aS^3}}{2}\right\}^2.
\end{align*}

\el

\bp Set $B_\varepsilon=\frac{\int_{\R^3} |\nabla u_{\varepsilon}|^2 \mathrm{d} x}{\left(\int_{\R^3} |u_{\varepsilon}|^6 \mathrm{d} x\right)^{\frac{1}{3}}}$ and $C_\varepsilon=\frac{\| u_{\varepsilon}\|_d^2}{\left(\int_{\R^3} |u_{\varepsilon}|^6 \mathrm{d} x\right)^{\frac{1}{3}}}$. By Lemma $2.1$, we get
$B_\varepsilon=S+O(\varepsilon^{\frac{1}{2}})$ and $C_\varepsilon=aS+O(\varepsilon^{\frac{1}{2}})$.
The direct calculation implies that
\begin{align}\label{2.7}
&\sup_{t \ge 0}\left[\frac{1}{2} t^2 \|u_\varepsilon\|_d^2+\frac{b}{4}t^4\left(\int_{\R^3}|\nabla u_\varepsilon|^2 \mathrm{d} x\right)^2-\frac{1}{6}t^6\int_{\R^3}|u_\varepsilon|^6\mathrm{d}x\right] \notag\\
&=\frac{1}{3}\left[\frac{b}{2}B_\varepsilon^2C_\varepsilon+\frac{1}{2}C_\varepsilon\sqrt{b^2B_\varepsilon^4+4C_\varepsilon}\right]+\frac{b}{12}
\left[\frac{b}{2}B_\varepsilon^3+\frac{1}{2}B_\varepsilon\sqrt{b^2B_\varepsilon^4+4C_\varepsilon}\right]^2\notag\\
&=\hat{c}+O(\varepsilon^{\frac{1}{2}}).
\end{align}
By Lemma $2.1$, there exists $\varepsilon_0>0$ such that for $\varepsilon \in (0,\varepsilon_0)$,
\begin{align*}
\int_{\R^3}|\nabla u_\varepsilon|^2 \mathrm{d} x \le \frac{3}{2} K_1, \ \ \ \int_{\R^3}| u_\varepsilon|^6 \mathrm{d} x \ge \frac{1}{2} K_2.
\end{align*}
Then by $I_d(t u_\varepsilon) \le \frac{t^2}{2} \|u_\varepsilon\|_d^2 + \frac{bt^4}{4}\left(\int_{\R^3}|\nabla u_\varepsilon|^2 \mathrm{d} x\right)^2 -\frac{t^6}{6}\int_{\mathbb{R}^{3}}|u_\varepsilon|^6\mathrm{d} x$, we can choose a small $t_1>0$ and a large $t_2>0$ such that
$\sup_{t \in [0, t_1]\cup[t_2,+\infty)}I_d(tu_\varepsilon)<\hat{c}$
independent of $\varepsilon \in (0,\varepsilon_0)$. By $(f_2)$, we have $\lim_{u \rightarrow +\infty}\frac{F(u)}{u^4}=+\infty$.
Then for any $l>0$, there exists $r_l>0$ such that $F(u) \ge l|u|^4$ for $|u| \ge r_l$.
Since $u_\varepsilon(x) \ge \varepsilon^{-\frac{1}{4}}$ for $|x| \le \varepsilon^{\frac{1}{2}} \le r$,
we derive $F(tu_\varepsilon) \ge lt_1^4 u_\varepsilon^4 \ge lt_1^4\varepsilon^{-1}$ for $t \in [t_1, t_2]$ and $|x| \le \varepsilon^{\frac{1}{2}} $ with $\varepsilon$ small enough. Together with $F(tu_\varepsilon) \ge 0$, we get
\begin{align}\label{2.8}
\inf_{t \in [t_1,t_2]}\int_{\R^3}F(tu_\varepsilon)\mathrm{d} x \ge & \inf_{t \in [t_1,t_2]}\int_{|x| \le \varepsilon^{\frac{1}{2}} }F(tu_\varepsilon)\mathrm{d} x \ge lt_1^4\varepsilon^{-1}\int_{|x| \le \varepsilon^{\frac{1}{2}}}1\mathrm{d} x = \eta_0 l \varepsilon^{\frac{1}{2}},
\end{align}
where $\eta_0>0$ is a constant. Combining (\ref{2.7})-(\ref{2.8}), we derive for $\varepsilon >0$ small,
\begin{align*}
&\sup_{t \in [t_1, t_2]}I_d(tu_\varepsilon)\\
 &\le   \sup_{t \ge 0}\left[\frac{1}{2} t^2 \|u_\varepsilon\|_d^2+\frac{b}{4}t^4\left(\int_{\R^3}|\nabla u_\varepsilon|^2 \mathrm{d} x\right)^2-\frac{1}{6}t^6\int_{\R^3}|u_\varepsilon|^6\mathrm{d}x\right] -\eta_0 l \varepsilon^{\frac{1}{2}}\\
&\le \hat{c} +C \varepsilon^{\frac{1}{2}}-\eta_0 l \varepsilon^{\frac{1}{2}}.
\end{align*}
By choosing $l$ large enough, we get $\sup_{t \in [t_1, t_2]}I_d(tu_\varepsilon)<\hat{c}$. So $\sup_{t \ge 0}I_d(tu_\varepsilon)<\hat{c}$ for $\varepsilon>0$ small. By the definition of $m_d$, we have $m_d \le \sup_{t \ge 0}I_d(tu_\varepsilon)$. Thus, we get $m_d<\hat{c}$.

\ep

By the argument of Proposition $2.6$ in \cite{HLP} and Lemma $2.2$, we know problem (\ref{2.3}) admits a positive ground state solution $u_d$. So
$m_d$ is attained by $u_d \in M_d$. Let $m_{r,d}=\{u \in H_r^1(\R^3) \setminus \{0\}:
(I_d'(u),u)=0\}$ with $M_{r,d}=\{u \in H_r^1(\R^3) \setminus \{0\}:
(I_d'(u),u)=0\}$. We can also derive $m_{r,d}$ is attained by $u_{r,d} \in M_{r,d}$. The proof is omitted here.

\bl\lab{Lemma 2.1} $m_{d}=m_{r,d}$.
\el

\bp Denote $H_r^1(\R^3)=H_r$ for simplicity. Set $c_d = \inf_{\gamma \in \Gamma} \max_{t \in [0,1]}I_d(\gamma(t))$,
where $\Gamma = \{\gamma \in C ([0,1],H): \gamma(0)=0, I_d(\gamma(1)) <0\}$. By Lemmas $4.1$-$4.2$ in \cite{Willem}, we can get the following equivalent characterization of $c_d$
\begin{align*}
m_d=c_d=\inf_{u \in H \setminus \{0\}} \sup_{t \ge 0} I_d(t u).
\end{align*}
Set the Pohozaev manifold
\begin{align*}
P=\left\{u \in H \setminus \{0\}: \frac{a}{2}\int_{\R^3}|\nabla u|^2 \mathrm{d}x+ \frac{b}{2}\left(\int_{\R^3}|\nabla u|^2 \mathrm{d}x\right)^2-3 \int_{\R^3}G(u) \mathrm{d}x=0\right\},
\end{align*}
where $G(u)=F(u)+\frac{1}{6}|u|^6-\frac{d}{2}|u|^2$. Then by the proof of Theorem $1.3$ in \cite{LG1}, we have
$c_d=\inf_{u \in P}I_d(u)$. Since the proof is standard, we omit it here.
So
\begin{align*}
m_d=c_d =\inf_{u \in P}I_d(u).
\end{align*}
Similarly, we also have $m_{r,d}=\inf_{u \in P_r}I_d(u)$, where
\begin{align*}
P_r=\left\{u \in H_r \setminus \{0\}: \frac{a}{2}\int_{\R^3}|\nabla u|^2 \mathrm{d}x+ \frac{b}{2}\left(\int_{\R^3}|\nabla u|^2 \mathrm{d}x\right)^2-3 \int_{\R^3}G(u) \mathrm{d}x=0\right\}.
\end{align*}
Recall that $\inf_{u \in P}I_d(u)=m_d=I_d(u_d)$.
Let $u_d^{\ast}$ the Schwarz spherical rearrangement of $u_d$. Then $u_d^{\ast}\in H_r$, $\int_{\R^3}|\nabla u_d^{\ast}|^2 \mathrm{d} x \le \int_{\R^3}|\nabla u_d|^2 \mathrm{d} x$ and $\int_{\R^3} G(u_d^{\ast}) \mathrm{d}x =\int_{\R^3} G(u_d) \mathrm{d}x$.
So
\begin{align}\label{2.4}
\frac{a}{2}\int_{\R^3}|\nabla u_d^{\ast}|^2 \mathrm{d}x+ \frac{b}{2}\left(\int_{\R^3}|\nabla u_d^{\ast}|^2 \mathrm{d}x\right)^2 \le 3 \int_{\R^3}G(u_d^{\ast}) \mathrm{d}x.
\end{align}
Since $\int_{\R^3}|\nabla u_d^{\ast}(\frac{.}{t})|^2 \mathrm{d}x=t \int_{\R^3}|\nabla u_d^{\ast}|^2 \mathrm{d}x$ and $\int_{\R^3}G(u_d^{\ast}(\frac{.}{t})) \mathrm{d}x=t^3\int_{\R^3}G(u_d^{\ast}) \mathrm{d}x$, we obtain that there exists $\hat{t}>0$ such that $u_d^{\ast} (\frac{.}{\hat{t}}) \in P_r$, that is,
\begin{align}\label{2.5}
\frac{a\hat{t}}{2}\int_{\R^3}|\nabla u_d^{\ast}|^2 \mathrm{d}x+ \frac{b\hat{t}^2}{2}\left(\int_{\R^3}|\nabla u_d^{\ast}|^2 \mathrm{d}x\right)^2 = 3 \hat{t}^3\int_{\R^3}G(u_d^{\ast}) \mathrm{d}x.
\end{align}
We claim $\hat{t} \le 1$. Otherwise, we have $\hat{t}>1$. Then by (\ref{2.4})-(\ref{2.5}),
\begin{align*}
\frac{a}{2}\int_{\R^3}|\nabla u_d^{\ast}|^2 \mathrm{d}x+ \frac{b}{2}\left(\int_{\R^3}|\nabla u_d^{\ast}|^2 \mathrm{d}x\right)^2 >& 3 \hat{t} \int_{\R^3}G(u_d^{\ast}) \mathrm{d}x \\
>&\frac{a}{2}\int_{\R^3}|\nabla u_d^{\ast}|^2 \mathrm{d}x+ \frac{b}{2}\left(\int_{\R^3}|\nabla u_d^{\ast}|^2 \mathrm{d}x\right)^2,
\end{align*}
a contradiction. So $\hat{t} \le 1$. By $u_d^{\ast}(\frac{.}{\hat{t}}) \in P_r$, $\int_{\R^3}|\nabla u_d^{\ast}(\frac{.}{\hat{t}})|^2 \mathrm{d}x=\hat{t} \int_{\R^3}|\nabla u_d^{\ast}|^2 \mathrm{d}x$ with $\hat{t} \le 1$ and $u_d \in P$,
\begin{align*}
m_{r,d}=\inf_{u \in P_r}I_d(u) \le & I_d\left(u_d^{\ast}(\frac{.}{\hat{t}})\right)\\
=&\frac{a}{3}\int_{\R^3}|\nabla u_d^{\ast}(\frac{.}{\hat{t}})|^2 \mathrm{d}x+ \frac{b}{12}\left(\int_{\R^3}|\nabla u_d^{\ast}(\frac{.}{\hat{t}})|^2 \mathrm{d}x\right)^2\\
\le & \frac{a}{3}\int_{\R^3}|\nabla u_d^{\ast}|^2 \mathrm{d}x+ \frac{b}{12}\left(\int_{\R^3}|\nabla u_d^{\ast}|^2 \mathrm{d}x\right)^2\\
\le & \frac{a}{3}\int_{\R^3}|\nabla u_d|^2 \mathrm{d}x+ \frac{b}{12}\left(\int_{\R^3}|\nabla u_d|^2 \mathrm{d}x\right)^2
=  I_d\left(u_d\right) = m_d.
\end{align*}
On the other hand, we have
\begin{align*}
m_d=\inf_{u \in P}I_d(u) \le \inf_{u \in P_r}I_d(u)=m_{r,d}.
\end{align*}
So $m_d=m_{r,d}$.
\ep

\bl\lab{Lemma 2.2} $m_\varepsilon \ge m_{V_0}$, $m_d \ge m_{V_0}$ if $d \ge V_0$.
\el

\bp For simplicity, we only proof $m_\varepsilon \ge m_{V_0}$. By the definition of $m_\varepsilon$, for any $\delta>0$, there exists $u_\delta \in M_\varepsilon$ such that
$m_\varepsilon>I_\varepsilon(u_\delta)-\delta$. By $(f_2)$, we know there exists a unique $t_\delta>0$ such that $t_\delta u_\delta \in
M_{V_0}$. By $t_\delta u_\delta \in M_{V_0}$, $u_\delta \in M_\varepsilon$ and $V(\varepsilon x) \ge V_0$,
\begin{align}\label{2-1}
&t_\delta^2\|u_\delta\|_{V_0}^2+bt_\delta^4\left(\int_{\R^3}|\nabla u_\delta|^2 \mathrm{d} x\right)^2=\int_{\mathbb{R}^{3}}f(t_\delta u_\delta) (t_\delta u_\delta)\mathrm{d} x+
t_\delta^6\int_{\mathbb{R}^{3}}|u_\delta|^6 \mathrm{d} x,\notag\\
&\|u_\delta\|_{V_0}^2+b\left(\int_{\R^3}|\nabla u_\delta|^2 \mathrm{d} x\right)^2 \le \int_{\mathbb{R}^{3}}f(u_\delta) u_\delta\mathrm{d} x+
\int_{\mathbb{R}^{3}}|u_\delta|^6 \mathrm{d} x.
\end{align}
If $t_\delta>1$, by $(f_2)$, we have $\int_{\mathbb{R}^{3}}f(t_\delta u_\delta) (t_\delta u_\delta)\mathrm{d} x > t_\delta^4 \int_{\mathbb{R}^{3}}f(u_\delta)u_\delta \mathrm{d} x$. Then
\begin{align*}
 t_\delta^4 \left(\int_{\mathbb{R}^{3}}f(u_\delta)u_\delta \mathrm{d} x+
\int_{\mathbb{R}^{3}}|u_\delta|^6 \mathrm{d} x\right)
< & t_\delta^4 \left[\|u_\delta\|_{V_0}^2+b\left(\int_{\R^3}|\nabla u_\delta|^2 \mathrm{d} x\right)^2\right]\\
\le & t_\delta^4 \left(\int_{\mathbb{R}^{3}}f(u_\delta)u_\delta \mathrm{d} x+
\int_{\mathbb{R}^{3}}|u_\delta|^6 \mathrm{d} x\right),
\end{align*}
a contradiction. So $t_\delta \le 1$, from which we get
\begin{align*}
\int_{\mathbb{R}^{3}}\left(\frac{1}{4}f(u_\delta)u_\delta-F(u_\delta)\right)
\mathrm{d}x \ge \int_{\mathbb{R}^{3}}\left(\frac{1}{4}f(t_\delta u_\delta)(t_\delta u_\delta)-F(t_\delta u_\delta)\right)
\mathrm{d}x.
\end{align*}
Thus,
\begin{align*}
m_\varepsilon+\delta& >  I_\varepsilon(u_\delta) -\frac{1}{4} \left(I_\varepsilon'(u_\delta),u_\delta\right)\\
& =
\frac{1}{4}\|u_\delta\|_{\varepsilon}^2+\int_{\mathbb{R}^{3}}\left(\frac{1}{4}f(u_\delta)u_\delta-F(u_\delta)\right)
\mathrm{d}x+\frac{1}{12}\int_{\mathbb{R}^{3}}|u_\delta|^6\mathrm{d}x\notag\\
& \ge
\frac{1}{4}t_\delta^2\|u_\delta\|_{V_0}^2+\int_{\mathbb{R}^{3}}\left(\frac{1}{4}f(t_\delta u_\delta)(t_\delta u_\delta)-F(t_\delta u_\delta)\right)
\mathrm{d}x+\frac{t_\delta^6}{12}\int_{\mathbb{R}^{3}}|u_\delta|^6\mathrm{d}x\notag\\
& = I_{V_0}(t_\delta u_\delta) -\frac{1}{4} \left(I'_{V_0}(t_\delta u_\delta),t_\delta u_\delta\right) =  I_{V_0}(t_\delta u_\delta) \ge m_{V_0}.
\end{align*}
Let $\delta \rightarrow 0$, we get $m_\varepsilon \ge m_{V_0}$.
\ep

\bl\lab{Lemma 2.3} Let $\{u_n\} \subset H$ be a sequence such that $I_d(u_n) \rightarrow c$, $I_d'(u_n) \rightarrow 0$ and $u_n \rightharpoonup u \ne 0$ weakly in $H$. Define the functionals $\hat{I}_d$, $\tilde{I}_d$ on $H$ by
\begin{align*}
\hat{I}_d(u)=&\frac{1}{2} \|u\|_d^2+\frac{bA}{4}\int_{\R^3}|\nabla u|^2 \mathrm{d} x-\int_{\R^3}F(u)\mathrm{d}x-\frac{1}{6}\int_{\R^3}|u|^6\mathrm{d}x,\\
\tilde{I}_d(u) =&\frac{1}{2} \|u\|_d^2+\frac{bA}{2}\int_{\R^3}|\nabla u|^2 \mathrm{d} x-\int_{\R^3}F(u)\mathrm{d}x-\frac{1}{6}\int_{\R^3}|u|^6\mathrm{d}x,
\end{align*}
where $A \ge \lim_{n \rightarrow \infty}\int_{\R^3}|\nabla u_{n}|^2 \mathrm{d} x$. Then $\tilde{I}_d'(u)=0$
and $c \ge \hat{I}_d(u) \ge m_d$.
\el

\bp By $I_d(u_n) \rightarrow c$ and $I_d'(u_n) \rightarrow 0$, we have
\begin{align*}
\hat{I}_d(u_{n})=c+o_n(1), \ \ \ \ \tilde{I}_d'(u_{n})=o_n(1).
\end{align*}
Then by $u_{n} \rightharpoonup u \neq 0$ weakly in $H$, we get $\tilde{I}'_d(u)=0$. Remark that
\begin{align*}
c+o_n(1)&=\hat{I}_d(u_{n}) -\frac{1}{4} \left(\tilde{I}_d'(u_{n}), u_n\right)\\
&=\frac{1}{4}\|u_{n}\|_d^2+\int_{\mathbb{R}^{3}}\left(\frac{1}{4}f(u_{n})u_{n}-F(u_{n})\right)
\mathrm{d}x+\frac{1}{12}\int_{\mathbb{R}^{3}}|u_{n}|^6\mathrm{d}x.
\end{align*}
Then by Fatou's Lemma,
\begin{align*}
c \ge & \frac{1}{4}\|u\|_d^2+\int_{\mathbb{R}^{3}}\left(\frac{1}{4}f(u)u-F(u)\right)
\mathrm{d}x+\frac{1}{12}\int_{\mathbb{R}^{3}}|u|^6\mathrm{d}x\\
=& \hat{I}_d(u) -\frac{1}{4} \left(\tilde{I}_d'(u), u\right)=\hat{I}_d(u).
\end{align*}
Since $u_n \rightharpoonup u \ne 0$ weakly in $H$, we have $A \ge \int_{\R^3}|\nabla u|^2 \mathrm{d} x$.
By $(f_2)$, we know there exists a unique $t>0$ such that $tu \in
M_d$, that is,
\begin{equation*}
t^2\|u\|_d^2+bt^4\left(\int_{\R^3}|\nabla u|^2 \mathrm{d} x\right)^2=\int_{\mathbb{R}^{3}}f(tu) (tu)\mathrm{d} x+
t^6\int_{\mathbb{R}^{3}}|u|^6 \mathrm{d} x.
\end{equation*}
By $(\tilde{I}'_d(u),u)=0$ and $A \ge \int_{\R^3}|\nabla u|^2 \mathrm{d} x$,
\begin{equation*}
\|u\|_d^2+b\left(\int_{\R^3}|\nabla u|^2 \mathrm{d} x\right)^2 \le \int_{\mathbb{R}^{3}}f(u)u \mathrm{d} x+
\int_{\mathbb{R}^{3}}|u|^6 \mathrm{d} x.
\end{equation*}
Similar to the argument of (\ref{2-1}), we have $t \le 1$. By $tu \in M_d$ and $t \le 1$,
\begin{align*}
\hat{I}_d(u) & =\hat{I}_d(u) -\frac{1}{4} \left(\tilde{I}'_d(u),u\right)\\
& \ge
\frac{1}{4}t^2\|u\|_d^2+\int_{\mathbb{R}^{3}}\left(\frac{1}{4}f(tu)(tu)-F(tu)\right)
\mathrm{d}x+\frac{t^6}{12}\int_{\mathbb{R}^{3}}|u|^6\mathrm{d}x\notag\\
& = I_d(tu) -\frac{1}{4} \left(I'_d(tu),tu\right) =  I_d(tu) \ge m_d.
\end{align*}

\ep

Now we introduce the barycenter function, which is crucial for proving multiplicity of solutions of (\ref{2.1}).
Consider the map $\Phi: H \rightarrow H$ defined by
\begin{equation*}
\Phi(u)(x):= \frac{1}{|B_1(x)|}\int_{B_1(x)}|u(y)| \mathrm{d}y, \ \ \forall \ x \in \R^3,
\end{equation*}
where $|B_1(x)|$ is the Lebesgue measure of $B_1(x)$. Set
\begin{equation*}
\hat{u}(x)=\left[\Phi(u)(x)-\frac{1}{2} \max_{x \in
\mathbb{R}^3}\Phi(u)(x)\right]^{+}.
\end{equation*}
Then we define the
barycenter $\beta: H \setminus \{0\}\rightarrow \mathbb{R}^3$ by
$\beta(u)=\frac{1}{\|\hat{u}\|_1}\int_{\mathbb{R}^3}x\hat{u}(x)\mathrm{d}x$.
From \cite{B-W,CP}, we know the map $\beta$ is continuous in $H \setminus \{0\}$ and
satisfies the following properties.

\vskip0.123in

\bl\lab{Lemma 2.4}
$\beta(u)=0$ if $u$ is radial; $\beta(tu)=\beta(u)$ for $t \ne 0$; $\beta(u(x-z))=\beta(u)+z$ for $z \in \mathbb{R}^3$.
\el

\s{Multiplicity of solutions of (\ref{1.1})}

\renewcommand{\theequation}{3.\arabic{equation}}

In this section, we study multiplicity of solutions of (\ref{1.1}). Since (\ref{1.1}) is equivalent to (\ref{2.1}), we consider (\ref{2.1}) instead.
For $l>0$, denote $C_l(x^i)$ the hypercube $\Pi_{j=1}^3(x_j^i-l,x_j^i-l)$ centered at $x^i=(x_1^i,x_2^i,x_3^i)$, $i=1$, $2$, \ldots $k$. Denote $\overline{C_l(x^i)}$ and $\partial C_l(x^i)$ the closure and the boundary of $C_l(x^i)$, respectively. By $(V_1)$-$(V_2)$, we can choose $l$, $L>0$
such that $\overline{C_l(x^i)}$, $i=1$, $2$, \ldots $k$ are disjoint, $V(x)>V(x^i)$ for $x \in  \overline{C_l(x^i)} \setminus x^i$, $i=1$, $2$, \ldots $k$
and $\cup_{i=1}^k \overline{C_l(x^i)} \subset \Pi_{j=1}^3(-L,L)$.

Let $C_{\frac{l}{\varepsilon}}(\frac{x^i}{\varepsilon})=C_{\frac{l}{\varepsilon}}^i$ and
\begin{align*}
&N_\varepsilon^i=\left\{u \in M_\varepsilon: u \ge 0, \ \beta(u) \in C_{\frac{l}{\varepsilon}}^i\right\},\\
& \partial N_\varepsilon^i=\left\{u \in M_\varepsilon: u \ge 0, \ \beta(u) \in \partial C_{\frac{l}{\varepsilon}}^i\right\},
\end{align*}
$i=1$, $2$, \ldots $k$. It can be readily verified that $N_\varepsilon^i$ and $\partial N_\varepsilon^i$ are non-empty sets for
$i=1$, $2$, \ldots $k$. Let
\begin{align*}
\gamma_\varepsilon^i= \inf_{u \in N_\varepsilon^i}I_\varepsilon (u), \ \ \ \ \tilde{\gamma}_\varepsilon^i= \inf_{u \in \partial N_\varepsilon^i}I_\varepsilon (u),
\end{align*}
$i=1$, $2$, \ldots $k$.
\vskip0.13in

\bl\lab{Lemma 3.1} For any $\delta \in (0,m_{V_0})$, there exists $\varepsilon_\delta>0$ such that for $\varepsilon \in (0,\varepsilon_\delta)$,
\begin{align*}
m_\varepsilon \le \gamma_\varepsilon^i < m_{V_0}+\delta,
\end{align*}
$i=1$, $2$, \ldots $k$. In particular, $N_\varepsilon^i$, $i=1$, $2$, \ldots $k$ are non-empty sets.
 \el

\bp Fix $i=1$, $2$, \ldots $k$. It is obvious that $\gamma_\varepsilon^i \ge m_\varepsilon$. By Lemma $2.3$, we know $m_{V_0}$ can be attained by a radial function $u_{r,V_0}$.
Let $\varepsilon \in (0,1)$. Define $\varphi_\varepsilon \in C_0^1(\R^3)$ such that $\varphi_\varepsilon(x)=1$ for $|x| < \frac{1}{\sqrt{\varepsilon}}-1$, $\varphi_\varepsilon(x)=0$ for $|x| > \frac{1}{\sqrt{\varepsilon}}$, $0 \le \varphi_\varepsilon \le 1$ and $|\nabla \varphi_\varepsilon| \le 2$. Let $w_\varepsilon^i(x)=u_{r,V_0}\left(x-\frac{x^i}{\varepsilon}\right)\varphi_\varepsilon\left(x-\frac{x^i}{\varepsilon}\right)$.
Set $g_\varepsilon(t)=I_\varepsilon(t w_\varepsilon^i)$, where $t >0$. We have
\begin{equation*}
g_\varepsilon'(t)=t\|w_\varepsilon^i\|_\varepsilon^2 + b t^3 \left(\int_{\R^3}|\nabla w_\varepsilon^i|^2 \mathrm{d} x\right)^2-\int_{\mathbb{R}^3}f(t
w_\varepsilon^i)w_\varepsilon^i\mathrm{d}
x-t^5\int_{\mathbb{R}^3}|w_\varepsilon^i|^6\mathrm{d} x.
\end{equation*}
Then by $(f_2)$, we know $g_\varepsilon(t)$ admits a unique critical point $t_\varepsilon^i >0$
corresponding to its maximum, that is, $g_\varepsilon(t_\varepsilon^i)=\sup_{t \ge 0}I_\varepsilon(t w_\varepsilon^i)$ and $g_\varepsilon'(t_\varepsilon^i)=0$. So $t_\varepsilon^i w_\varepsilon^i \in M_\varepsilon$. By Lemma $2.6$, we get $\beta\left(t_\varepsilon^i w_\varepsilon^i\right)=\beta\left(w_\varepsilon^i\right)=\frac{x^i}{\varepsilon}+\beta\left(u_{r,V_0}\varphi_\varepsilon\right)$.
Since $\lim_{\varepsilon \rightarrow 0}\|u_{r,V_0}\varphi_\varepsilon-u_{r,V_0}\|_{V_0}=0$,
by the continuity of $\beta$, we get $\beta\left(u_{r,V_0}\varphi_\varepsilon\right) \rightarrow \beta\left(u_{r,V_0}\right)=0$ as $\varepsilon \rightarrow 0$,
in view of $u_{r,V_0}$ is radial. So $\beta\left(t_\varepsilon^i w_\varepsilon^i\right)=\frac{x^i}{\varepsilon}+O(\varepsilon)$, where $O(\varepsilon) \rightarrow 0$ as $\varepsilon \rightarrow 0$, from which we derive $\beta\left(t_\varepsilon^i w_\varepsilon^i\right) \in C_{\frac{l}{\varepsilon}}^i$ for $\varepsilon>0$ small.
Recall that $t_\varepsilon^i w_\varepsilon^i \in M_\varepsilon$. Then $t_\varepsilon^i w_\varepsilon^i \in N_\varepsilon^i$, that is, $N_\varepsilon^i$ is non-empty.

By $t_\varepsilon^i w_\varepsilon^i \in N_\varepsilon^i$, we have $\gamma_\varepsilon^i \le I_\varepsilon(t_\varepsilon^i w_\varepsilon^i)=\sup_{t \ge 0}I_\varepsilon(t w_\varepsilon^i)$. Recall that $\lim_{\varepsilon \rightarrow 0}\|u_{r,V_0}\varphi_\varepsilon-u_{r,V_0}\|_{V_0}=0$. Then
\begin{align}\label{3.1}
\lim_{\varepsilon \rightarrow 0}\int_{\R^3}F(w_\varepsilon^i)\mathrm{d}x=&\lim_{\varepsilon \rightarrow 0}\int_{\R^3}F(u_{r,V_0}\varphi_\varepsilon)\mathrm{d}x= \int_{\R^3}F(u_{r,V_0})\mathrm{d}x.
\end{align}
Similarly, we have
\begin{align}\label{3.2}
&\lim_{\varepsilon \rightarrow 0}\int_{\R^3}f(w_\varepsilon^i)w_\varepsilon^i\mathrm{d}x= \int_{\R^3}f(u_{r,V_0})u_{r,V_0}\mathrm{d}x,\notag\\
&\lim_{\varepsilon \rightarrow 0}\int_{\R^3}|\nabla w_\varepsilon^i|^2\mathrm{d}x=\int_{\R^3}|\nabla u_{r,V_0}|^2\mathrm{d}x,\ \ \
\lim_{\varepsilon \rightarrow 0}\int_{\R^3}|w_\varepsilon^i|^6\mathrm{d}x= \int_{\R^3}|u_{r,V_0}|^6\mathrm{d}x.
\end{align}
By the Lebesgue dominated convergence theorem and $V(x^i)=V_0$, we also have
\begin{align}\label{3.3}
\lim_{\varepsilon \rightarrow 0}\int_{\R^3}V(\varepsilon x) |w_\varepsilon^i|^2\mathrm{d}x=&\lim_{\varepsilon \rightarrow 0}\int_{\R^3}V(\varepsilon x+x^i) |u_{r,V_0}\varphi_\varepsilon|^2\mathrm{d}x=\int_{\R^3}V_0 |u_{r,V_0}|^2\mathrm{d}x.
\end{align}
We claim $\lim_{\varepsilon \rightarrow 0}t_\varepsilon^i=1$. In fact, by $t_\varepsilon^i w_\varepsilon^i \in N_\varepsilon^i$,
\begin{align}\label{3.4}
(t_\varepsilon^i)^2\|w_\varepsilon^i\|_{\varepsilon}^2+b(t_\varepsilon^i)^4\left(\int_{\R^3}|\nabla w_\varepsilon^i|^2 \mathrm{d}x\right)^2=
\int_{\R^3}\left(f(t_\varepsilon^i w_\varepsilon^i)t_\varepsilon^i w_\varepsilon^i+(t_\varepsilon^i)^6|w_\varepsilon^i|^6\right)\mathrm{d}x.
\end{align}
Then
$\frac{1}{(t_\varepsilon^i)^4}\|w_\varepsilon^i\|_{\varepsilon}^2+\frac{b}{(t_\varepsilon^i)^2}\left(\int_{\R^3}|\nabla w_\varepsilon^i|^2 \mathrm{d}x\right)^2 \ge \int_{\R^3}|w_\varepsilon^i|^6\mathrm{d}x$.
If $t_\varepsilon^i \rightarrow +\infty$, by (\ref{3.2})-(\ref{3.3}), we get a contradiction. So $t_\varepsilon^i$ is bounded. By $(f_1)$,
we derive for $\eta=\frac{V_0}{2}$, there exists $C_\eta=C_{\frac{V_0}{2}}>0$ such that
\begin{align*}
(t_\varepsilon^i)^2\|w_\varepsilon^i\|_{\varepsilon}^2 \le \frac{(t_\varepsilon^i)^2}{2} \int_{\R^3}V(\varepsilon x)|w_\varepsilon^i|^2 \mathrm{d}x +C_{\frac{V_0}{2}}(t_\varepsilon^i)^6 \int_{\R^3}|w_\varepsilon^i|^6 \mathrm{d}x,
\end{align*}
in view of $V(\varepsilon x) \ge V_0$ and (\ref{3.4}).
Then $\frac{1}{2} \|w_\varepsilon^i\|_{\varepsilon}^2 \le C_{\frac{V_0}{2}}(t_\varepsilon^i)^4 \int_{\R^3}|w_\varepsilon^i|^6 \mathrm{d}x$,
that is, $(t_\varepsilon^i)^4 \ge \frac{\|w_\varepsilon^i\|_{\varepsilon}^2 }{2 C_{\frac{V_0}{2}} \int_{\R^3}|w_\varepsilon^i|^6 \mathrm{d}x}$.
So by (\ref{3.2})-(\ref{3.3}), we can assume $t_\varepsilon^i \rightarrow t^i > 0$. Let $\varepsilon \rightarrow 0$ in (\ref{3.4}),
we have
\begin{align*}
&(t^i)^2\|u_{r,V_0}\|_{V_0}^2+b(t^i)^4\left(\int_{\R^3}|\nabla u_{r,V_0}|^2 \mathrm{d}x\right)^2\\
&=
\int_{\R^3}\left(f(t^i u_{r,V_0})(t^i u_{r,V_0})+(t^i)^6|u_{r,V_0}|^6\right)\mathrm{d}x,
\end{align*}
that is, $t^i u_{r,V_0} \in M_{V_0}$. Note that $u_{r,V_0} \in M_{V_0}$ and there exists a unique $t>0$ satisfying $t u_{r,V_0} \in M_{V_0}$. Then
 $t^i=1$. Combining (\ref{3.1})-(\ref{3.3}) and $t_\varepsilon^i \rightarrow 1$,
\begin{align*}
\gamma_\varepsilon^i \le  I_\varepsilon(t_\varepsilon^i w_\varepsilon^i)= I_\varepsilon(w_\varepsilon^i)+O(\varepsilon)=I_{V_0}(u_{r,V_0})+O(\varepsilon)=m_{V_0}+O(\varepsilon).
\end{align*}
So for any $\delta \in (0,m_{V_0})$, there exists $\varepsilon_\delta>0$ such that $\gamma_\varepsilon^i < m_{V_0}+\delta$ for $\varepsilon \in (0,\varepsilon_\delta)$.
\ep

\bl\lab{Lemma 3.2} There exist $\eta$, $\varepsilon_\eta>0$ such that $\widetilde{\gamma}_\varepsilon^i > m_{V_0}+ \eta$ for $\varepsilon \in (0,\varepsilon_\eta)$, $i=1$, $2$, \ldots $k$.

\el

\bp Fix $i=1$, $2$, \ldots $k$. Assume by contradiction that there exists a sequence $\{\varepsilon_n\}$ such that $\varepsilon_n \rightarrow 0$ and
$\widetilde{\gamma}_{\varepsilon_n}^i \rightarrow \tilde{c} \le m_{V_0}$. Then there exists $\{u_n\} \subset \partial N_{\varepsilon_n}^i$ such that
$I_{\varepsilon_n}(u_n) \rightarrow \tilde{c} \le m_{V_0}$. Since $\{u_n\} \subset \partial N_{\varepsilon_n}^i$, we have $\beta(u_n) \in \partial C_{\frac{l}{\varepsilon_n}}^i$ and
\begin{align*}
\|u_n\|_{V_0}^2+b\left(\int_{\R^3}|\nabla u_n|^2 \mathrm{d} x\right)^2 \le &\|u_n\|_{\varepsilon_n}^2+b\left(\int_{\R^3}|\nabla u_n|^2 \mathrm{d} x\right)^2\notag\\
=& \int_{\R^3}f(u_n) u_n\mathrm{d}x+ \int_{\R^3}|u_n|^6\mathrm{d}x.
\end{align*}
By $(f_2)$, there exists $t_n>0$ such that $t_n u_n \in M_{V_0}$, that is,
\begin{align*}
t_n^2\|u_n\|_{V_0}^2+b t_n^4 \left(\int_{\R^3}|\nabla u_n|^2 \mathrm{d} x\right)^2 = \int_{\R^3}f(t_n u_n)t_n u_n\mathrm{d}x+t_n^6 \int_{\R^3}|u_n|^6\mathrm{d}x.
\end{align*}
Similar to the argument of (\ref{2-1}), we have $t_n \le 1$. Then
\begin{align*}
m_{V_0} \ge & I_{\varepsilon_n}(u_n) - \frac{1}{4} \left(I_{\varepsilon_n}'(u_n),u_n\right)+o_n(1)\\
\ge &I_{V_0}(u_n)-\frac{1}{4} \left(I_{V_0}'(u_n),u_n\right) +o_n(1)\\
\ge & I_{V_0}(t_n u_n)-\frac{1}{4} \left(I_{V_0}'(t_nu_n),t_nu_n\right) +o_n(1)\\
=& I_{V_0}(t_n u_n) +o_n(1) \ge m_{V_0}+o_n(1),
\end{align*}
from which we get $t_n \rightarrow 1$ and
\begin{align}\label{3.5}
\int_{\R^3}V(\varepsilon_n x)|u_n|^2\mathrm{d}x=\int_{\R^3}V_0|u_n|^2\mathrm{d}x+o_n(1).
\end{align}
Moreover, set $\bar{u}_n=t_nu_n$, we get
\begin{align}\label{3-0}
I_{V_0}(\bar{u}_n)=m_{V_0}+o_n(1), \ \ \ \left(I_{V_0}'(\bar{u}_n),\bar{u}_n\right)=0.
\end{align}
Let $G_{V_0}(u)=(I_{V_0}'(u),u)$. By the Ekeland's variational principle, there exist $\{v_n\} \subset M_{V_0}$ and $\mu_n \in \R^1$ such that
$\|v_n-\bar{u}_n\|_{V_0} = o_n(1)$, $I_{V_0}(v_n)=m_{V_0}+o_n(1)$ and $I_{V_0}'(v_n)-\mu_n G_{V_0}'(v_n)=o_n(1)$.
Then $\mu_n (G_{V_0}'(v_n),v_n)=o_n(1)$. Since
\begin{equation*}
(G_{V_0}'(v_n),v_n)=2\|v_n\|_{V_0}^2 + 4 b  \left(\int_{\R^3}|\nabla v_n|^2 \mathrm{d} x\right)^2- \int_{\R^3}( f(v_n)v_n +f'(v_n)v_n^2+6 |v_n|^6)\mathrm{d} x,
\end{equation*}
by $f'(v_n)v_n^2-3f(v_n)v_n \ge 0$, we have
\begin{align}\label{3.6}
(G_{V_0}'(v_n),v_n)=&(G_{V_0}'(v_n),v_n) -4(I_{V_0}'(v_n),v_n) \notag\\
\le &2 \|v_n\|_{V_0}^2+4 b \left(\int_{\R^3}|\nabla v_n|^2 \mathrm{d} x\right)^2-4  \int_{\mathbb{R}^3}  f(
v_n) v_n \mathrm{d} x -6 \int_{\mathbb{R}^3}|v_n|^6\mathrm{d}x\notag\\
&-4(I_{V_0}'(v_n),v_n)\notag\\
=&-2 \|v_n\|_{V_0}^2-2 \int_{\mathbb{R}^3}|v_n|^6\mathrm{d}x<0,
\end{align}
from which we get $\lim_{n \rightarrow \infty}(G_{V_0}'(v_n),v_n) \le 0$. If
$\lim_{n \rightarrow \infty}(G_{V_0}'(v_n),v_n) = 0$, by (\ref{3.6}), we have $\|v_n\|_{V_0} \rightarrow 0$,
a contradiction with $I_{V_0}(v_{n}) \rightarrow m_{V_0}>0$. So $\lim_{n \rightarrow \infty}(G'(v_n),v_n) < 0$.
By $\mu_n (G'(v_n),v_n)=o_n(1)$, we get $I_{V_0}'(v_n) \rightarrow 0$. Thus,
\begin{align}\label{3.7}
I_{V_0}(v_n)=m_{V_0}+o_n(1), \ \ \ \ I_{V_0}'(v_n)=o_n(1).
\end{align}
By (\ref{3.7}), we can derive $\|v_n\|_{V_0}$ is bounded. Assume $\lim_{n \rightarrow \infty} \int_{\R^3}|\nabla v_n|^2 \mathrm{d} x$ exist.
By the Lions Lemma, we know $\int_{\R^3}|v_n|^t\mathrm{d}x \rightarrow 0$ for any $t \in (2,6)$, or there exists $y_{n} \in \mathbb{R}^{3}$ such that $v_{n}(.+y_{n}) \rightharpoonup w \neq 0$ weakly in $H$.
If $\int_{\R^3}|v_n|^t\mathrm{d}x \rightarrow 0$ for any $t \in (2,6)$, by $(f_1)$, we have $\int_{\R^3}F(v_n)\mathrm{d}x \rightarrow 0$ and $\int_{\R^3}f(v_n)v_n\mathrm{d}x \rightarrow 0$.
Then by (\ref{3.7}), we derive
\begin{align}\label{3.9}
&m_{V_0} + o_n(1)= \frac{1}{2}\|v_n\|_{V_0}^2+\frac{b}{4}\left(\int_{\R^3}|\nabla v_n|^2 \mathrm{d} x\right)^2-\frac{1}{6}
\int_{\mathbb{R}^{3}}|v_n|^6\mathrm{d}x,\notag\\
& \|v_n\|_{V_0}^2 +b\left(\int_{\R^3}|\nabla v_n|^2 \mathrm{d} x\right)^2- \int_{\mathbb{R}^{3}}|v_n|^6\mathrm{d}x = o_n(1).
\end{align}
So
\begin{align}\label{3.10}
&m_{V_0} \ge \frac{a}{3}\int_{\R^3}|\nabla v_n|^2 \mathrm{d} x+\frac{b}{12}\left(\int_{\R^3}|\nabla v_n|^2 \mathrm{d} x\right)^2+ o_n(1),\notag\\
& a\int_{\R^3}|\nabla v_n|^2 \mathrm{d} x +b\left(\int_{\R^3}|\nabla v_n|^2 \mathrm{d} x\right)^2 \le \int_{\mathbb{R}^{3}}|v_n|^6\mathrm{d}x + o_n(1).
\end{align}
If $\lim_{n \rightarrow \infty} \int_{\R^3}|\nabla v_n|^2 \mathrm{d} x=0$, by the Sobolev embedding theorem, we have
$\lim_{n \rightarrow \infty} \int_{\R^3}|v_n|^6 \mathrm{d} x=0$. Then by (\ref{3.9}), we get $\|v_n\|_{V_0} \rightarrow 0$, a contradiction with $m_{V_0}>0$.
So $\lim_{n \rightarrow \infty} \int_{\R^3}|\nabla v_n|^2 \mathrm{d} x>0$.
By the second inequality in (\ref{3.10}) and
$\int_{\mathbb{R}^{3}}|v_n|^6\mathrm{d}x \le \frac{\left(\int_{\R^3}|\nabla v_n|^2 \mathrm{d} x\right)^3}{S^3}$,
we derive $\lim_{n \rightarrow \infty} \int_{\R^3}|\nabla v_n|^2 \mathrm{d} x \ge \frac{bS^3+\sqrt{(bS^3)^2+4aS^3}}{2}$.
Then by the first inequality in (\ref{3.10}), we have $m_{V_0} \ge \hat{c}$, a contradiction with $m_{V_0}<\hat{c}$. Here $\hat{c}$ is defined in Lemma $2.2$. Thus, we derive $v_{n}(.+y_{n}) \rightharpoonup w \neq 0$ weakly in $H$.
Together with (\ref{3.7}), we have
\begin{align*}
I_{V_0}(v_{n}(.+y_{n}))=m_{V_0}+o_n(1), \ \ \ I_{V_0}'(v_{n}(.+y_{n}))=o_n(1).
\end{align*}
Let $A= \lim_{n \rightarrow \infty}\int_{\R^3}|\nabla v_{n}(.+y_{n})|^2 \mathrm{d} x$.
Define the functionals $\hat{I}_{V_0}$, $\tilde{I}_{V_0}$ on $H$ by
\begin{align*}
\hat{I}_{V_0}(u)=&\frac{1}{2} \|u\|_{V_0}^2+\frac{bA}{4}\int_{\R^3}|\nabla u|^2 \mathrm{d} x-\int_{\R^3}F(u)\mathrm{d}x-\frac{1}{6}\int_{\R^3}|u|^6\mathrm{d}x,\\
\tilde{I}_{V_0}(u) =&\frac{1}{2} \|u\|_{V_0}^2+\frac{bA}{2}\int_{\R^3}|\nabla u|^2 \mathrm{d} x-\int_{\R^3}F(u)\mathrm{d}x-\frac{1}{6}\int_{\R^3}|u|^6\mathrm{d}x.
\end{align*}
Then
\begin{align*}
\hat{I}_{V_0}(v_{n}(.+y_{n}))=m_{V_0}+o_n(1), \ \ \ \tilde{I}_{V_0}'(v_{n}(.+y_{n}))=o_n(1).
\end{align*}
By $v_{n}(.+y_{n}) \rightharpoonup w$ weakly in $H$ and Lemma $2.5$, we have $\tilde{I}'_{V_0}(w)=0$ and $\hat{I}_{V_0}(w) \ge m_{V_0}$.
Observe that
\begin{align*}
m_{V_0}+o_n(1) =  & \hat{I}_{V_0}(v_n(.+y_{n})) -\frac{1}{4}(\tilde{I}_{V_0}'(v_n(.+y_{n})),v_n(.+y_{n})) \\
=& \frac{1}{4} \|v_n(.+y_{n})\|_{V_0}^2+\frac{1}{12} \int_{\R^3}|v_n(.+y_{n})|^6 \mathrm{d} x\\
&+ \int_{\R^3}\left(\frac{1}{4}f(v_n(.+y_{n}))v_n(.+y_{n})-F(v_n(.+y_{n}))\right)\mathrm{d} x.
\end{align*}
By Fatou's Lemma,
\begin{align*}
m_{V_0}\ge  &  \frac{1}{4} \|w\|_{V_0}^2+\frac{1}{12} \int_{\R^3}|w|^6 \mathrm{d} x+ \int_{\R^3}\left(\frac{1}{4}f(w)w-F(w)\right)\mathrm{d} x\\
=& \hat{I}_{V_0}(w) -\frac{1}{4}(\tilde{I}_{V_0}'(w),w)=\hat{I}_{V_0}(w) \ge m_{V_0}.
\end{align*}
So $v_{n}(.+y_{n}) \rightarrow w$ in $H$. By the continuity of $\beta$ and Lemma $2.6$, we have $\beta (v_{n})-y_n =\beta (v_{n}(.+y_{n})) \rightarrow \beta(w)$.
Then by $\|v_n-\bar{u}_n\|_{V_0} \rightarrow 0$, we get $\beta (\bar{u}_{n})-y_n \rightarrow \beta(w)$. So $\mathrm{dist} \left(\varepsilon_n y_n, \partial C_{l}(x^i)\right) \rightarrow 0$ in view of $\beta(\bar{u}_n)=\beta(u_n) \in \partial C_{\frac{l}{\varepsilon_n}}^i$. Assume $\varepsilon_n y_n \rightarrow y_0 \in \partial C_{l}(x^i)$. Then $V(y_0) > V_0 =\inf_{x \in \R^3}V(x)$. By $t_n \rightarrow 1$ and (\ref{3.5}),
\begin{align*}
\int_{\R^3}V(\varepsilon_n x)|\bar{u}_n|^2\mathrm{d}x=\int_{\R^3}V_0|\bar{u}_n|^2\mathrm{d}x+o_n(1).
\end{align*}
Recall that $\|v_n-\bar{u}_n\|_{V_0} \rightarrow 0$ and $v_{n}(.+y_{n}) \rightarrow w$ in $H$. Then we have
\begin{align*}
\lim_{n \rightarrow \infty}\int_{\R^3}V_0|\bar{u}_n|^2\mathrm{d}x=&\lim_{n \rightarrow \infty}\int_{\R^3}V_0|v_n|^2\mathrm{d}x\\
=&\lim_{n \rightarrow \infty}\int_{\R^3}V_0|v_n(.+y_{n})|^2\mathrm{d}x=\int_{\R^3}V_0|w|^2\mathrm{d}x.
\end{align*}
We also have $\bar{u}_{n}(.+y_{n}) \rightarrow w$ in $H$. Then by Fatou's Lemma,
\begin{align*}
\lim_{n \rightarrow \infty}\int_{\R^3}V(\varepsilon_n x)|\bar{u}_n|^2\mathrm{d}x = & \lim_{n \rightarrow \infty} \int_{\R^3}V(\varepsilon_n x+ \varepsilon_n y_n)|\bar{u}_{n}(.+y_{n})|^2\mathrm{d}x\\
\ge & \int_{\R^3}V(y_0)|w|^2\mathrm{d}x,
\end{align*}
in view of $\varepsilon_n y_n \rightarrow y_0$. Thus, we obtain $\int_{\R^3}V(y_0)|w|^2\mathrm{d}x \le \int_{\R^3}V_0|w|^2\mathrm{d}x$, a contradiction with $V(y_0) > V_0$.

\ep

Let
\begin{align*}
G_\varepsilon(u)=&(I_\varepsilon'(u),u)=\|u\|_\varepsilon^2 +b \left(\int_{\R^3}|\nabla u|^2 \mathrm{d} x\right)^2 - \int_{R^3}\left( f(u)u+ |u|^6\right)\mathrm{d} x.
\end{align*}
Then for any $\varphi \in H_\varepsilon$,
\begin{align*}
(G_\varepsilon'(u),\varphi)=&2(u,\varphi)_\varepsilon + 4 b \int_{\R^3}|\nabla u|^2 \int_{\R^3}\nabla u \nabla \varphi \mathrm{d} x\\
&- \int_{\R^3}( f(u) +f'(u)u+6 u^5 )\varphi\mathrm{d} x.
\end{align*}
Following the idea of \cite{NT}, we get the following result.

\vskip0.1in
\bl\lab{Lemma 3.3} For any $u \in N_\varepsilon^i$, $i=1$, $2$, \ldots $k$, there exist a positive constant $\sigma$ and a differential function
$s(w)>0$ with $w \in H_\varepsilon$ and $\|w\|_\varepsilon < \sigma$, satisfying
\begin{itemize}
\item [($i$)] $s(0)=1 \ \mathrm{and} \
s(w)(u+w) \in N_\varepsilon^i \  \mathrm{for} \ \mathrm{any} \ \|w\|_\varepsilon <\sigma$;
\item [($ii$)] $(s'(0), \varphi)=\frac{-(G_\varepsilon'(u),\varphi)}{(G_\varepsilon'(u),u)}$ for any $\varphi \in H_\varepsilon$.
\end{itemize}
\el

\bp Fix $i=1$, $2$, \ldots $k$. For any $u \in N_\varepsilon^i$, define
$T: H_\varepsilon \times \mathbb{R} \rightarrow \mathbb{R}$ by
\begin{align*}
T(w,s)=& s^2 \|u+w\|_\varepsilon^2+b s^4 \left(\int_{\R^3}|\nabla (u+w)|^2 \mathrm{d} x\right)^2-  \int_{\mathbb{R}^3}  f(
su + sw) (su+sw) \mathrm{d} x \\
&- s^6 \int_{\mathbb{R}^3}| u+ w|^6\mathrm{d}x.
\end{align*}
Since $u \in N_\varepsilon^i$, we have $T(0,1)=0$. Then by
$f'(u)u^2-3f(u)u \ge 0$,
\begin{align*}
T_s(0,1)=&2 \|u\|_\varepsilon^2+4 b \left(\int_{\R^3}|\nabla u|^2 \mathrm{d} x\right)^2-  \int_{\mathbb{R}^3}  f(
u) u \mathrm{d} x -  \int_{\mathbb{R}^3}  f'(u) u^2 \mathrm{d} x\\
& -6 \int_{\mathbb{R}^3}|u|^6\mathrm{d}x\\
\le &2 \|u\|_\varepsilon^2+4 b \left(\int_{\R^3}|\nabla u|^2 \mathrm{d} x\right)^2-4  \int_{\mathbb{R}^3}  f(
u) u \mathrm{d} x -6 \int_{\mathbb{R}^3}|u|^6\mathrm{d}x,
\end{align*}
from which we derive
$T_{s}(0,1)=T_{s}(0,1)-4(I_\varepsilon'(u),u)<0$.
By the implicit function theorem at the point $(0,1)$, we obtain that there exist a positive constant $\sigma$ and a differential function
$s(w)>0$ with $w \in H_\varepsilon$ and $\|w\|_\varepsilon < \sigma$, satisfying $s(w)(u+w) \in M_\varepsilon$ for $\|w\|_\varepsilon <\sigma$.
Since $u \in N_\varepsilon^i$, we have $\beta(u) \in C_{\frac{l}{\varepsilon}}^i$. Then by the continuity of functions $\beta$ and $s$,
we get $\beta(s(w)(u+w)) \in C_{\frac{l}{\varepsilon}}^i$ for $\sigma>0$ small. So $s(w)(u+w) \in N_\varepsilon^i$.
Moreover, since $s(w)$ is a differential function, by the direct calculation, we get $(ii)$.

\ep

\bl\lab{Lemma 3.4} Fix $i=1$, $2$, \ldots $k$. Then there exists a sequence $\{u_{n}\} \subset N_\varepsilon^i$ such
that $I_\varepsilon(u_{n}) \rightarrow \gamma_\varepsilon^i$ and $I_\varepsilon'(u_{n}) \rightarrow 0$.
\el

\bp By the definition of $\gamma_\varepsilon^i$, there is $\{u_n\}
\subset N_\varepsilon^i$ such that $I_\varepsilon(u_n) \rightarrow \gamma_\varepsilon^i$.  Then
we derive $\|u_n\|_\varepsilon$ is bounded.
By the Ekeland's variational principle,
\begin{align*}
& I_\varepsilon(u_{n}) \leq  \gamma_\varepsilon^i + \frac{1}{n}, \ \
&I_\varepsilon (v) \ge I_\varepsilon (u_{n}) - \frac{1}{n}
\|u_{n} - v\|_\varepsilon, \ \ \  \forall \ v \in N_\varepsilon^i.
\end{align*}
By Lemma $3.3$, there
exist $\sigma_n \downarrow 0$ and $s_{n} (w)$ satisfying
\begin{equation*}
s_{n} (w) (u_{n} + w) \in N_\varepsilon^i, \  \  \forall\    w \in H_\varepsilon, \  \  \|w\|_\varepsilon < \sigma_{n}.
\end{equation*}
Let $w = t \phi$, where $\phi \in H_\varepsilon$ and $t > 0$ small. Then
\begin{align*}
&\frac{1}{n} \left[t s_{n} (t \phi) \|\phi\|_\varepsilon + \left|s_{n}(t \phi)
-1\right| \|u_{n}\|_\varepsilon\right]\\
& \ge \frac{1}{n} \left\|u_{n} -
s_{n}(t \phi) (u_{n} + t
\phi)\right\|_\varepsilon\\
& \ge I_\varepsilon(u_{n}) - I_\varepsilon\left(s_{n} (t \phi) (u_{n} + t \phi)\right)\\
&= \left[I_\varepsilon(u_{n}) - I_\varepsilon(u_{n} + t \phi)\right] + \left[I_\varepsilon(u_{n} + t
\phi) - I_\varepsilon(s_{n} (t
\phi) (u_{n} + t \phi))\right]\\
&= \left[I_\varepsilon(u_{n}) - I_\varepsilon(u_{n} + t \phi)\right]\\
&\quad + (1-s_{n} (t
\phi))\left(I_\varepsilon'(u_{n} + t
\phi +\theta_n (s_{n} (t
\phi) (u_{n} + t \phi))), u_{n} + t \phi\right),
\end{align*}
where $\theta_n \in (0,1)$. Dividing by $t$ and let $t \rightarrow 0$, we derive
\begin{align*}
\frac{1}{n}\left[\left|(s_{n}'(0), \phi)\right|\|u_{n}\|_\varepsilon + \|\phi\|_\varepsilon\right]
 \ge -(I_\varepsilon'(u_{n}),\phi),
\end{align*}
in view of $(I_\varepsilon'(u_n),u_n)=0$. By Lemma $3.3$, we have
$(s_n'(0), \phi)=\frac{-(G_\varepsilon'(u_n),\phi)}{(G_\varepsilon'(u_n),u_n)}$.
We claim $\lim_{n \rightarrow \infty}(G_\varepsilon'(u_n),u_n) < 0$.
In fact, by $f'(u_n)u_n^2-3f(u_n)u_n \ge 0$,
\begin{align*}
&(G_\varepsilon'(u_n),u_n)\\
&=(G_\varepsilon'(u_n),u_n) -4(I_\varepsilon'(u_n),u_n) \\
& \le 2 \|u_n\|_\varepsilon^2+4 b \left(\int_{\R^3}|\nabla u_n|^2 \mathrm{d} x\right)^2-4  \int_{\mathbb{R}^3}  f(
u_n) u_n \mathrm{d} x -6 \int_{\mathbb{R}^3}|u_n|^6\mathrm{d}x\\
&\quad -4(I_\varepsilon'(u_n),u_n)\\
&=-2 \|u_n\|_\varepsilon^2-2 \int_{\mathbb{R}^3}|u_n|^6\mathrm{d}x<0.
\end{align*}
So $\lim_{n \rightarrow \infty}(G_\varepsilon'(u_n),u_n) \le 0$. Moreover, if
$\lim_{n \rightarrow \infty}(G_\varepsilon'(u_n),u_n) = 0$, then $u_n \rightarrow 0$ in $H_\varepsilon$,
a contradiction with $I_\varepsilon(u_{n}) \rightarrow \gamma_\varepsilon^i>0$. By $\lim_{n \rightarrow \infty}(G'(u_n),u_n)<0$, $\|u_n\|_\varepsilon$ is bounded
and $(s_n'(0), \phi)=\frac{-(G_\varepsilon'(u_n),\phi)}{(G_\varepsilon'(u_n),u_n)}$, we derive there exists $M>0$ such that $|(s_n'(0), \phi)| \le M$ for any $n$.
So $\frac{1}{n}\left(M \|u_{n}\|_\varepsilon + \|\phi\|_\varepsilon\right)
 \ge -(I_\varepsilon'(u_{n}),\phi)$. Let $n \rightarrow \infty$, we get $I_\varepsilon'(u_n) \rightarrow 0$.

\ep

\vskip0.1in

From Lemmas $3.1$-$3.2$, we know there exists $\hat{\varepsilon}>0$ such that for $\varepsilon \in (0,\hat{\varepsilon})$,
\begin{align*}
m_\varepsilon \le \gamma_\varepsilon^i <\tilde{\gamma}_\varepsilon^i.
\end{align*}
By Lemmas $2.4$ and $3.1$, we have $\gamma_\varepsilon^i \rightarrow m_{V_0}$ as $\varepsilon \rightarrow 0$. Together with Lemma $2.2$, we derive there exists $\tilde{\varepsilon} \in (0,\hat{\varepsilon})$ such that
\begin{align*}
\gamma_\varepsilon^i<\min\{ \hat{c}, 2 m_{V_0}\}
\end{align*}
for $\varepsilon \in (0,\tilde{\varepsilon})$.

\vskip0.13in

\bl\lab{Lemma 3.5} Fix $i=1$, $2$, \ldots $k$. Let $\{u_{n}\} \subset N_\varepsilon^i$ be a sequence such
that $I_\varepsilon(u_{n}) \rightarrow \gamma_\varepsilon^i$ and $I_\varepsilon'(u_{n}) \rightarrow 0$.
Then $\{u_{n}\}$ converges strongly in $H_\varepsilon$ up to a subsequence
for $\varepsilon \in (0,\tilde{\varepsilon})$. \el

\bp Obviously, we have $\|u_n\|_\varepsilon$ is bounded. Assume $u_{n} \rightharpoonup u$
weakly in $H_\varepsilon$. We claim $u \ne 0$. Otherwise, we have $u_{n} \rightharpoonup 0$ weakly in $H_\varepsilon$.
By the definition of $V_\infty$, for any $\delta>0$, there exists $R_\delta>0$ such that
$V(\varepsilon x) \ge V_\infty - \delta$ for $|x| \ge R_\delta$. So by $u_{n} \rightharpoonup 0$ weakly in $H_\varepsilon$,
\begin{align*}
\int_{\R^3}V(\varepsilon x)|u_n|^2 \mathrm{d} x =& \int_{|x| \ge R_\delta}V(\varepsilon x)|u_n|^2 \mathrm{d} x+\int_{|x| \le R_\delta}V(\varepsilon x)|u_n|^2 \mathrm{d} x\\
\ge & \int_{|x| \ge R_\delta}V_\infty|u_n|^2 \mathrm{d} x - \delta \int_{|x| \ge R_\delta}|u_n|^2 \mathrm{d} x+ o_n(1)\\
\ge & \int_{\R^3}V_\infty|u_n|^2 \mathrm{d} x - \delta \int_{\R^3}|u_n|^2 \mathrm{d} x+ o_n(1).
\end{align*}
Let $\delta \rightarrow 0$, we have
\begin{align}\label{3-1}
\int_{\R^3}V(\varepsilon x)|u_n|^2 \mathrm{d} x \ge \int_{\R^3}V_\infty|u_n|^2 \mathrm{d} x+o_n(1).
\end{align}
Then by $u_n \in M_\varepsilon$, we get
\begin{align}\label{3--1}
&\lim_{n \rightarrow \infty} \int_{\mathbb{R}^{3}}f(u_n) u_n\mathrm{d} x+
\lim_{n \rightarrow \infty} \int_{\mathbb{R}^{3}}|u_n|^6 \mathrm{d} x\notag\\
& \ge \lim_{n \rightarrow \infty} \left[\|u_n\|_{V_\infty}^2+b\left(\int_{\R^3}|\nabla u_n|^2 \mathrm{d} x\right)^2\right].
\end{align}
By $(f_2)$, there exists $t_n>0$ such that $t_n u_n \in M_{V_\infty}$, that is,
\begin{align}\label{3-2}
t_n^2\|u_n\|_{V_\infty}^2+bt_n^4\left(\int_{\R^3}|\nabla u_n|^2 \mathrm{d} x\right)^2= &\int_{\mathbb{R}^{3}}f(t_nu_n) (t_nu_n)\mathrm{d} x+
t_n^6\int_{\mathbb{R}^{3}}|u_n|^6 \mathrm{d} x\notag\\
\ge & t_n^6\int_{\mathbb{R}^{3}}|u_n|^6 \mathrm{d} x.
\end{align}
By $u_n \in M_\varepsilon$ and $(f_1)$, we derive for $\eta=\frac{V_0}{2}$, there exists $C_\eta=C_{\frac{V_0}{2}}$ such that
\begin{align*}
\|u_n\|_\varepsilon^2 \le & \frac{V_0}{2}  \int_{\mathbb{R}^{3}}|u_n|^2\mathrm{d} x+\left(C_{\frac{V_0}{2}}+1\right)
\int_{\mathbb{R}^{3}}|u_n|^6 \mathrm{d} x\\
\le & \frac{1}{2}  \|u_n\|_\varepsilon^2+\left(C_{\frac{V_0}{2}}+1\right)
\int_{\mathbb{R}^{3}}|u_n|^6 \mathrm{d} x.
\end{align*}
If $\lim_{n \rightarrow \infty}\int_{\mathbb{R}^{3}}|u_n|^6 \mathrm{d} x=0$, then $\|u_n\|_\varepsilon^2 \rightarrow 0$, a contradiction with
$\gamma_\varepsilon^i>0$. So $\lim_{n \rightarrow \infty}\int_{\mathbb{R}^{3}}|u_n|^6 \mathrm{d} x>0$. Together with (\ref{3-2}),
we derive $t_n$ is bounded. Assume $t_n \rightarrow t_0$. We claim $t_0 \le 1$. If $t_0 >1$, without loss of generality, we may assume $t_n>1$ for any $n \in N$. So $\int_{\mathbb{R}^{3}}f(t_nu_n) (t_nu_n)\mathrm{d} x > t_n^4 \int_{\mathbb{R}^{3}}f(u_n)u_n \mathrm{d} x$. Together with (\ref{3-2}), we have
\begin{align}\label{3-3}
& t_0^4 \lim_{n \rightarrow \infty}\int_{\mathbb{R}^{3}}f(u_n)u_n \mathrm{d} x+
 t_0^6 \lim_{n \rightarrow \infty} \int_{\mathbb{R}^{3}}|u_n|^6 \mathrm{d} x\notag\\
& \le  t_0^4 \lim_{n \rightarrow \infty} \left[\|u_n\|_{V_\infty}^2+b\left(\int_{\R^3}|\nabla u_n|^2 \mathrm{d} x\right)^2\right].
\end{align}
Combining (\ref{3--1}), (\ref{3-3}), $t_0>1$ and $\lim_{n \rightarrow \infty}\int_{\mathbb{R}^{3}}|u_n|^6 \mathrm{d} x>0$,
\begin{align*}
&\lim_{n \rightarrow \infty} \int_{\mathbb{R}^{3}}f(u_n) u_n\mathrm{d} x+
\lim_{n \rightarrow \infty} \int_{\mathbb{R}^{3}}|u_n|^6 \mathrm{d} x\\
& \ge \lim_{n \rightarrow \infty} \left[\|u_n\|_{V_\infty}^2+b\left(\int_{\R^3}|\nabla u_n|^2 \mathrm{d} x\right)^2\right] \\
& \ge \lim_{n \rightarrow \infty} \int_{\mathbb{R}^{3}}f(u_n) u_n\mathrm{d} x+t_0^2
\lim_{n \rightarrow \infty} \int_{\mathbb{R}^{3}}|u_n|^6 \mathrm{d} x\\
& >  \lim_{n \rightarrow \infty} \int_{\mathbb{R}^{3}}f(u_n) u_n\mathrm{d} x+
\lim_{n \rightarrow \infty} \int_{\mathbb{R}^{3}}|u_n|^6 \mathrm{d} x,
\end{align*}
a contradiction. By $u_n \in M_\varepsilon$, $t_n u_n \in M_{V_\infty}$ with $t_n \rightarrow t_0 \le 1$ and (\ref{3-1}),
\begin{align*}
I_\varepsilon(u_n)& =  I_\varepsilon(u_n) -\frac{1}{4} \left(I_\varepsilon'(u_n),u_n\right)\\
& =
\frac{1}{4}\|u_n\|_{\varepsilon}^2+\int_{\mathbb{R}^{3}}\left(\frac{1}{4}f(u_n)u_n-F(u_n)\right)
\mathrm{d}x+\frac{1}{12}\int_{\mathbb{R}^{3}}|u_n|^6\mathrm{d}x\notag\\
& \ge
\frac{1}{4}t_n^2\|u_n\|_{V_\infty}^2+\int_{\mathbb{R}^{3}}\left(\frac{1}{4}f(t_n u_n)(t_n u_n)-F(t_n u_n)\right)
\mathrm{d}x\\
& \quad +\frac{t_n^6}{12}\int_{\mathbb{R}^{3}}|u_n|^6\mathrm{d}x+o_n(1)\notag\\
& = I_{V_\infty}(t_n u_n) -\frac{1}{4} \left(I'_{V_\infty}(t_n u_n),t_n u_n\right) =  I_{V_\infty}(t_n u_n)+o_n(1).
\end{align*}
So
\begin{align}\label{3_4}
I_\varepsilon(u_n) \ge  I_{V_\infty}(t_n u_n)+o_n(1), \ \ \  \left(I_{V_\infty}'(t_n u_n),t_n u_n\right)=0.
\end{align}
Moreover, by $u_{n} \rightharpoonup 0$ weakly in $H_\varepsilon$ and $t_n \rightarrow t_0 \le 1$, we get $t_n u_n \rightharpoonup 0$ weakly in $H_\varepsilon$.
Similar to the argument of (\ref{3-0}) and (\ref{3.7}), we can derive from (\ref{3_4}) that there exists $\{\check{u}_n\} \subset H$ satisfying $\|\check{u}_n-t_n u_n\|_{V_\infty} = o_n(1)$, $I_{V_\infty}(t_n u_n)=  I_{V_\infty}(\check{u}_n)+o_n(1)$ and $I_{V_\infty}'(\check{u}_n)=o_n(1)$. So
\begin{align}\label{3-4}
\gamma_\varepsilon^i= I_\varepsilon(u_n) + o_n(1)\ge  I_{V_\infty}(\check{u}_n)+o_n(1), \ \ \ \ \ \  I_{V_\infty}'(\check{u}_n)=o_n(1).
\end{align}
Since the embedding $H_\varepsilon \hookrightarrow H$ is continuous, we have $t_n u_n \rightharpoonup 0$ weakly in $H$. Then by $\|\check{u}_n-t_n u_n\|_{V_\infty} = o_n(1)$,
we get $\check{u}_n \rightharpoonup 0$ weakly in $H$.
By (\ref{3_4}) and $(f_1)$, for $\eta=\frac{V_\infty}{2}$, there exists $C_\eta=C_{\frac{V_\infty}{2}}$ such that
\begin{align*}
\|t_n u_n\|_{V_\infty}^2 \le & \frac{V_\infty}{2}  \int_{\mathbb{R}^{3}}|t_n u_n|^2\mathrm{d} x+\left(C_{\frac{V_\infty}{2}}+1\right)
\int_{\mathbb{R}^{3}}|t_n u_n|^6 \mathrm{d} x\\
\le & \frac{1}{2}  \|t_n u_n\|_{V_\infty}^2+\left(C_{\frac{V_\infty}{2}}+1\right)
\int_{\mathbb{R}^{3}}|t_n u_n|^6 \mathrm{d} x.
\end{align*}
Since $S \left(\int_{\mathbb{R}^{3}}|t_n u_n|^6 \mathrm{d} x\right)^{\frac{1}{3}} \le \|t_n u_n\|_{V_\infty}^2$, we get
$\int_{\mathbb{R}^{3}}|t_n u_n|^6 \mathrm{d} x \ge \left[\frac{S}{2\left(C_{\frac{V_\infty}{2}}+1\right)}\right]^{\frac{3}{2}}$.
Then by $\|\check{u}_n-t_n u_n\|_{V_\infty} = o_n(1)$, we have
\begin{align}\label{3--17}
\lim_{n \rightarrow \infty}\int_{\mathbb{R}^{3}}|\check{u}_n|^6 \mathrm{d} x \ge \left[\frac{S}{2\left(C_{\frac{V_\infty}{2}}+1\right)}\right]^{\frac{3}{2}}>0.
\end{align}
The Lions Lemma implies that $\int_{\R^3}|\check{u}_n|^t\mathrm{d}x \rightarrow 0$ for any $t \in (2,6)$, or there exists $z_{n} \in \mathbb{R}^{3}$ with $|z_n| \rightarrow \infty$ such that $u_n^1=\check{u}_{n}(.+z_{n}) \rightharpoonup u^1 \neq 0$ weakly in $H$.
Thus,
if $\int_{\R^3}|\check{u}_n|^t\mathrm{d}x \rightarrow 0$ for any $t \in (2,6)$, similar to the argument of (\ref{3.9})-(\ref{3.10}),
we can derive from (\ref{3-4})-(\ref{3--17}) that $\gamma_\varepsilon^i \ge \hat{c}$, a contradiction.
So $u_{n}^1=\check{u}_n(.+z_{n}) \rightharpoonup u^1 \neq 0$ weakly in $H$ with $|z_n| \rightarrow \infty$. By (\ref{3-4}), we have
\begin{align*}
\gamma_\varepsilon^i \ge I_{V_\infty}(u^1_{n}) + o_n(1), \ \ \ \ \ \ I_{V_\infty}'(u^1_{n}) = o_n(1).
\end{align*}
Since $\beta(t_n u_n)=\beta(u_n) \in C_{\frac{l}{\varepsilon}}^i$, by $\|\check{u}_n-t_n u_n\|_{V_\infty} \rightarrow 0$,
we have $\beta(\check{u}_n) \in C_{\frac{l}{\varepsilon}}^i$ for $n$ large enough.
Then by $\beta(u^1_n)=\beta(\check{u}_n)-z_n$, we have
\begin{align*}
\frac{x^i-l}{\varepsilon}-z_n<\beta(u^1_n)<\frac{x^i+l}{\varepsilon}-z_n
\end{align*}
for $n$ large enough. Let $\hat{A}= \lim_{n \rightarrow \infty}\int_{\R^3}|\nabla u^1_n|^2 \mathrm{d} x$.
Define the functionals $\hat{I}_{V_\infty}$, $\tilde{I}_{V_\infty}$ on $H$ by
\begin{align*}
\hat{I}_{V_\infty}(u)=&\frac{1}{2} \|u\|_{V_\infty}^2+\frac{b\hat{A}}{4}\int_{\R^3}|\nabla u|^2 \mathrm{d} x-\int_{\R^3}F(u)\mathrm{d}x-\frac{1}{6}\int_{\R^3}|u|^6\mathrm{d}x,\\
\tilde{I}_{V_\infty}(u) =&\frac{1}{2} \|u\|_{V_\infty}^2+\frac{b\hat{A}}{2}\int_{\R^3}|\nabla u|^2 \mathrm{d} x-\int_{\R^3}F(u)\mathrm{d}x-\frac{1}{6}\int_{\R^3}|u|^6\mathrm{d}x.
\end{align*}
Then
\begin{align}\label{3-5}
|\beta(u^1_n)| \rightarrow \infty, \ \ \gamma_\varepsilon^i \ge \hat{I}_{V_\infty}(u^1_{n}) + o_n(1), \ \ \tilde{I}_{V_\infty}'(u^1_{n}) = o_n(1).
\end{align}
Since $u^1_n \rightharpoonup u^1 \neq 0$ weakly in $H$, by Lemma $2.5$, we have
\begin{align*}
\hat{I}_{V_\infty}(u^1) \ge m_{V_\infty}, \ \ \ \ \tilde{I}_{V_\infty}'(u^1)=0.
\end{align*}
Set $v_{n}= u^1_{n}-u^1$. By Lemma $1.3$ in \cite{CH},
\begin{align}\label{319}
\int_{\R^3} F(u^1_n)\mathrm{d}x- \int_{\R^3} F(u^1)\mathrm{d}x= \int_{\R^3} F(v_n)\mathrm{d}x+o_n(1).
\end{align}
The Brezis-Lieb Lemma in $\cite{Willem}$ implies that
\begin{align}\label{320}
\int_{\R^3}|u_n^1|^6\mathrm{d}x-\int_{\R^3}|u^1|^6\mathrm{d}x=\int_{\R^3}|v_n|^6\mathrm{d}x+o_n(1).
\end{align}
Together with (\ref{3-5}), we get
$\gamma_\varepsilon^i \ge \hat{I}_{V_\infty}(v_{n})+ \hat{I}_{V_\infty}(u^1) + o_n(1)$.
On the other hand, by Lemma $8.9$ in $\cite{Willem}$, we know for any $\varphi \in H$,
\begin{align}\label{321}
\left|\int_{\mathbb{R}^3}\left[(u_{n}^1)^5 -(u^1)^5- (v_n)^{5}
\right]\varphi \mathrm{d} x\right| = o_n(1) \|\varphi\|_{V_\infty}.
\end{align}
Similar to Lemma $8.1$ in $\cite{Willem}$, we obtain that for any $\varphi \in H$,
\begin{align}\label{322}
\left|\int_{\mathbb{R}^3}\left[f(u_{n}^1) -f(u^1)- f(v_n)
\right]\varphi \mathrm{d} x\right| = o_n(1) \|\varphi\|_{V_\infty}.
\end{align}
Together with $\tilde{I}_{V_\infty}'(u^1_{n}) \rightarrow 0$ and $\tilde{I}_{V_\infty}'(u^1)=0$, we get $\tilde{I}_{V_\infty}'(v_{n}) \rightarrow 0$.
Thus,
\begin{align}\label{3.12}
\gamma_\varepsilon^i \ge \hat{I}_{V_\infty}(v_{n})+ \hat{I}_{V_\infty}(u^1) + o_n(1), \ \ \ \ \tilde{I}_{V_\infty}'(v_{n}) = o_n(1).
\end{align}
If $v_n \rightarrow 0$ in $H$, that is, $u^1_n \rightarrow u^1$ in $H$, by the continuity of $\beta$, we have $\beta(u^1_n) \rightarrow \beta(u^1)$, a contradiction with $|\beta(u^1_n)| \rightarrow \infty$. So $v_n$ converges weakly(not strongly) to $0$ in $H$.
The Lions Lemma implies that $\int_{\R^3}|v_n|^t\mathrm{d}x \rightarrow 0$ for any $t \in (2,6)$, or there exists $z_{n}^1 \in \mathbb{R}^{3}$ with $|z_n^1| \rightarrow \infty$ such that $v_n^1=v_{n}(.+z_{n}^1) \rightharpoonup v^1 \neq 0$ weakly in $H$.
If $\int_{\R^3}|v_n|^t\mathrm{d}x \rightarrow 0$ for any $t \in (2,6)$, by (\ref{3.12}),
\begin{align}\label{324}
&\lim_{n \rightarrow \infty}\hat{I}_{V_\infty}(v_{n})= \frac{1}{2}\|v_n\|_{V_\infty}^2+\frac{b\hat{A}}{4}\int_{\R^3}|\nabla v_n|^2 \mathrm{d} x-\frac{1}{6}
\int_{\mathbb{R}^{3}}|v_n|^6\mathrm{d}x + o_n(1),\notag\\
& \|v_n\|_{V_\infty}^2 +b\hat{A}\int_{\R^3}|\nabla v_n|^2 \mathrm{d} x- \int_{\mathbb{R}^{3}}|v_n|^6\mathrm{d}x = o_n(1).
\end{align}
By $\hat{A}= \lim_{n \rightarrow \infty}\int_{\R^3}|\nabla u^1_n|^2 \mathrm{d} x$ and
$v_{n}= u^1_{n}-u^1 \rightharpoonup 0$ weakly in $H$, we derive $\hat{A}\ge \lim_{n \rightarrow \infty}\int_{\R^3}|\nabla v_n|^2 \mathrm{d} x$.
Thus,
\begin{align}\label{325}
&\lim_{n \rightarrow \infty}\hat{I}_{V_\infty}(v_{n}) \ge \frac{a}{3}\lim_{n \rightarrow \infty}\int_{\R^3}|\nabla v_n|^2 \mathrm{d} x+\frac{b}{12}\left(\lim_{n \rightarrow \infty}\int_{\R^3}|\nabla v_n|^2 \mathrm{d} x\right)^2,\notag\\
& a\lim_{n \rightarrow \infty}\int_{\R^3}|\nabla v_n|^2 \mathrm{d} x +b\left(\lim_{n \rightarrow \infty}\int_{\R^3}|\nabla v_n|^2 \mathrm{d} x\right)^2 \le \lim_{n \rightarrow \infty}\int_{\mathbb{R}^{3}}|v_n|^6\mathrm{d}x.
\end{align}
If $\lim_{n \rightarrow \infty} \int_{\R^3}|\nabla v_n|^2 \mathrm{d} x=0$, by the Sobolev embedding theorem, we have
$\lim_{n \rightarrow \infty} \int_{\R^3}|v_n|^6 \mathrm{d} x=0$. Then by (\ref{324}), we get $\|v_n\|_{V_\infty} \rightarrow 0$, a contradiction with 
$v_n$ converges weakly(not strongly) to $0$ in $H$. So $\lim_{n \rightarrow \infty} \int_{\R^3}|\nabla v_n|^2 \mathrm{d} x>0$.
By (\ref{325}) and $\int_{\mathbb{R}^{3}}|v_n|^6\mathrm{d}x \le \frac{\left(\int_{\R^3}|\nabla v_n|^2 \mathrm{d} x\right)^3}{S^3}$,
we derive $\lim_{n \rightarrow \infty} \int_{\R^3}|\nabla v_n|^2 \mathrm{d} x \ge \frac{bS^3+\sqrt{(bS^3)^2+4aS^3}}{2}$.
So $\lim_{n \rightarrow \infty}\hat{I}_{V_\infty}(v_{n}) \ge \hat{c}$.
Recall that $\hat{I}_{V_\infty}(u^1) \ge m_{V_\infty}$. Then by (\ref{3.12}), we get
$\gamma_\varepsilon^i > \hat{c}$, a contradiction.
So $v_n^1=v_{n}(.+z_{n}^1) \rightharpoonup v^1 \neq 0$ weakly in $H$ with $|z_n^1| \rightarrow \infty$.
By (\ref{3.12}), we have
\begin{align}\label{3.13}
\gamma_\varepsilon^i \ge \hat{I}_{V_\infty}(v_{n}^1)+ \hat{I}_{V_\infty}(u^1) + o_n(1), \ \ \ \ \tilde{I}_{V_\infty}'(v_{n}^1) = o_n(1).
\end{align}
Since $\hat{A}\ge \lim_{n \rightarrow \infty}\int_{\R^3}|\nabla v_n|^2 \mathrm{d} x=\lim_{n \rightarrow \infty}\int_{\R^3}|\nabla v_n^1|^2 \mathrm{d} x$ and  $v_n^1 \rightharpoonup v^1 \neq 0$ weakly in $H$, by (\ref{3.13}) and Lemma $2.5$, we derive $\lim_{n \rightarrow \infty} \hat{I}_{V_\infty}(v_{n}^1) \ge m_{V_\infty}$.
Together with $\hat{I}_{V_\infty}(u^1) \ge m_{V_\infty}$, we get $\gamma_\varepsilon^i \ge 2 m_{V_\infty} \ge 2 m_{V_0}$, a contradiction.
So $u_{n} \rightharpoonup u \ne 0$ weakly in $H_\varepsilon$.

Let $\tilde{A}= \lim_{n \rightarrow \infty}\int_{\R^3}|\nabla u_n|^2 \mathrm{d} x$. Define the functionals $\hat{I}_\varepsilon$, $\tilde{I}_\varepsilon$ on $H_\varepsilon$ by
\begin{align*}
\hat{I}_\varepsilon(u)=&\frac{1}{2} \|u\|_\varepsilon^2+\frac{b\tilde{A}}{4}\int_{\R^3}|\nabla u|^2 \mathrm{d} x-\int_{\R^3}F(u)\mathrm{d}x-\frac{1}{6}\int_{\R^3}|u|^6\mathrm{d}x,\\
\tilde{I}_\varepsilon(u) =&\frac{1}{2} \|u\|_\varepsilon^2+\frac{b\tilde{A}}{2}\int_{\R^3}|\nabla u|^2 \mathrm{d} x-\int_{\R^3}F(u)\mathrm{d}x-\frac{1}{6}\int_{\R^3}|u|^6\mathrm{d}x.
\end{align*}
By $I_\varepsilon(u_{n}) \rightarrow \gamma_\varepsilon^i$ and $I_\varepsilon'(u_{n}) \rightarrow 0$, we have $\hat{I}_\varepsilon(u_{n}) \rightarrow \gamma_\varepsilon^i$ and $\tilde{I}_\varepsilon'(u_{n}) \rightarrow 0$. Since $u_{n} \rightharpoonup u \ne 0$
weakly in $H_\varepsilon$, similar to Lemma $2.5$, we get $\hat{I}_\varepsilon(u) \ge m_\varepsilon$ and $\tilde{I}_\varepsilon'(u) = 0$. Let $\tilde{u}_n=u_n-u$. Similar to (\ref{3-5})-(\ref{3.12}), we can derive from $I_\varepsilon(u_{n}) \rightarrow \gamma_\varepsilon^i$ and $I_\varepsilon'(u_{n}) \rightarrow 0$ that
\begin{align}\label{3.14}
\gamma_\varepsilon^i = \hat{I}_\varepsilon(\tilde{u}_{n}) + \hat{I}_\varepsilon(u) + o_n(1), \ \ \ \ \ \ \tilde{I}_\varepsilon'(\tilde{u}_{n}) = o_n(1).
\end{align}
We claim $\tilde{u}_{n} \rightarrow 0$ in $H_\varepsilon$. Otherwise, $\tilde{u}_n$ converges weakly(not strongly) to $0$ in $H_\varepsilon$.
Define the functionals $\acute{I}_{V_\infty}$, $\grave{I}_{V_\infty}$ on $H$ by
\begin{align*}
\acute{I}_{V_\infty}(u)=&\frac{1}{2} \|u\|_{V_\infty}^2+\frac{b\tilde{A}}{4}\int_{\R^3}|\nabla u|^2 \mathrm{d} x-\int_{\R^3}F(u)\mathrm{d}x-\frac{1}{6}\int_{\R^3}|u|^6\mathrm{d}x,\\
\grave{I}_{V_\infty}(u) =&\frac{1}{2} \|u\|_{V_\infty}^2+\frac{b\tilde{A}}{2}\int_{\R^3}|\nabla u|^2 \mathrm{d} x-\int_{\R^3}F(u)\mathrm{d}x-\frac{1}{6}\int_{\R^3}|u|^6\mathrm{d}x.
\end{align*}
Similar to (\ref{3_4}), we can derive from (\ref{3.14}) that there exists $\hat{t}_n>0$ satisfying $\tilde{t}_n \rightarrow \tilde{t} \le 1$ and
\begin{align*}
 \hat{I}_\varepsilon(\tilde{u}_{n}) \ge  \acute{I}_{V_\infty}(\tilde{t}_n \tilde{u}_n)+o_n(1), \ \ \  \left(\grave{I}_{V_\infty}'(\tilde{t}_n \tilde{u}_n),\tilde{t}_n \tilde{u}_n\right)=0.
\end{align*}
Then similar to (\ref{3-0}) and (\ref{3.7}), we obtain that there exists $\{\acute{u}_n\} \subset H$ such that
$\|\acute{u}_n-\tilde{t}_n \tilde{u}_n\|_{V_\infty} = o_n(1)$, $\acute{I}_{V_\infty}(\tilde{t}_n \tilde{u}_n)=  \acute{I}_{V_\infty}(\acute{u}_n)+o_n(1)$ and $\grave{I}_{V_\infty}'(\acute{u}_n)=o_n(1)$. So $\acute{u}_n \rightharpoonup 0$ weakly in $H$ and
\begin{align}\label{3-6}
\gamma_\varepsilon^i= \hat{I}_\varepsilon(\tilde{u}_{n}) + \hat{I}_\varepsilon(u) + o_n(1) \ge  \acute{I}_{V_\infty}(\acute{u}_n)+ \hat{I}_\varepsilon(u)+o_n(1), \ \   \grave{I}_{V_\infty}'(\acute{u}_n)=o_n(1).
\end{align}
The Lions Lemma implies that $\int_{\R^3}|\acute{u}_n|^t\mathrm{d}x \rightarrow 0$ for any $t \in (2,6)$, or there exists $\acute{z}_{n} \in \mathbb{R}^{3}$ with $|\acute{z}_n| \rightarrow \infty$ such that $\acute{u}_{n}^1=\acute{u}_{n}(.+\acute{z}_{n}) \rightharpoonup \acute{u}^1 \neq 0$ weakly in $H$.
If $\int_{\R^3}|\acute{u}_n|^t\mathrm{d}x \rightarrow 0$ for any $t \in (2,6)$, similar to (\ref{324})-(\ref{325}),
we can derive from (\ref{3-6}) that $\lim_{n \rightarrow \infty} \acute{I}_{V_\infty}(\acute{u}_{n})  \ge \hat{c}$.
So $\gamma_\varepsilon^i > \hat{c}$, a contradiction.
Then $\acute{u}_{n}^1=\acute{u}_n(.+\acute{z}_{n}) \rightharpoonup \acute{u}^1 \neq 0$ weakly in $H$ with $|\acute{z}_n| \rightarrow \infty$.
By (\ref{3-6}),
\begin{align}\label{3.15}
\gamma_\varepsilon^i \ge  \acute{I}_{V_\infty}(\acute{u}^1_{n}) +\hat{I}_\varepsilon(u) + o_n(1), \ \ \ \ \ \ \grave{I}_{V_\infty}'(\acute{u}^1_{n}) = o_n(1).
\end{align}
Since $\tilde{A}= \lim_{n \rightarrow \infty}\int_{\R^3}|\nabla u_n|^2 \mathrm{d} x$ and
$\tilde{u}_n=u_{n}-u \rightharpoonup 0$ weakly in $H_\varepsilon$, we have $\tilde{A} \ge \lim_{n \rightarrow \infty}\int_{\R^3}|\nabla \tilde{u}_n|^2 \mathrm{d} x$. Together with $\tilde{t}_n \rightarrow \tilde{t} \le 1$, we get $\tilde{A} \ge \lim_{n \rightarrow \infty}\int_{\R^3}|\nabla (t_n \tilde{u}_n)|^2 \mathrm{d} x$. Then by $\|\acute{u}_n-\tilde{t}_n \tilde{u}_n\|_{V_\infty} \rightarrow 0$, we obtain that
\begin{align*}
\tilde{A} \ge \lim_{n \rightarrow \infty}\int_{\R^3}|\nabla \acute{u}_n|^2 \mathrm{d} x = \lim_{n \rightarrow \infty}\int_{\R^3}|\nabla \acute{u}_n^1|^2 \mathrm{d} x.
\end{align*}
Since $\acute{u}_{n}^1 \rightharpoonup \acute{u}^1 \neq 0$ weakly in $H$, by Lemma $2.5$, we derive from (\ref{3.15}) that $\lim_{n \rightarrow \infty}\acute{I}_{V_\infty}(\acute{u}^1_{n}) \ge m_{V_\infty} \ge m_{V_0}$. Together with $\hat{I}_\varepsilon(u) \ge m_\varepsilon \ge m_{V_0}$, we get $\gamma_\varepsilon^i \ge  2m_{V_0}$,
a contradiction. So $u_{n} \rightarrow u$ in $H_\varepsilon$.

\ep

\bl\lab{Lemma 3.6} Problem (\ref{2.1}) admits at least $k$ different solutions for $\varepsilon \in (0,\tilde{\varepsilon})$.
\el

\bp By Lemma $3.4$, for any fixed $i=1$, $2$, \ldots $k$, there exists $\{u_{n}\} \subset N_\varepsilon^i$ satisfying $I_\varepsilon(u_{n}) \rightarrow \gamma_\varepsilon^i$ and $I_\varepsilon'(u_{n}) \rightarrow 0$. Then by Lemma $3.5$,
$u_n^i \rightarrow u^i$ in $H_\varepsilon$. So $u^i \in N_\varepsilon^i \cup \partial N_\varepsilon^i$, $I_\varepsilon(u^i) = \gamma_\varepsilon^i$ and $I_\varepsilon'(u^i) = 0$. Since $\gamma_\varepsilon^i < \tilde{\gamma}_\varepsilon^i$, we get
$u^i \in N_\varepsilon^i$. By $C_{\frac{l}{\varepsilon}}^i$, $i=1$, $2$, \ldots $k$ are disjoint and $\beta(u^i) \in C_{\frac{l}{\varepsilon}}^i$, we obtain that $u^i$, $i=1$, $2$, \ldots $k$ are different. Obviously, $u^i$ is non-negative. The maximum principle implies that $u^i$ is positive. So (\ref{2.1}) admits at least $k$ different positive solutions.

\ep

\s{Concentration of solutions of (\ref{1.1})}

\renewcommand{\theequation}{4.\arabic{equation}}

Let $u_\varepsilon^i(x)$, $i=1$, $2$, \ldots $k$ be solutions of (\ref{2.1}). Then
$u_\varepsilon^i \in N_\varepsilon^i$, $I_\varepsilon(u_\varepsilon^i) = \gamma_\varepsilon^i$ and $I_\varepsilon'(u_\varepsilon^i) = 0$.
Now we study the concentration of $u_\varepsilon^i$ as $\varepsilon \rightarrow 0$.

\bl\lab{Lemma 4.1} Fix $i=1$, $2$, \ldots $k$. Then there exist $\varepsilon_0 \in (0,\tilde{\varepsilon})$, $\{x_\varepsilon^i\} \subset \R^3$,
$R_0>0$ and $\gamma_0>0$ such that $\int_{B_{R_0}(x_\varepsilon^i)}|u_\varepsilon^i|^2 \mathrm{d}x \ge \gamma_0$ for $\varepsilon \in (0,\varepsilon_0)$.
\el

\bp Otherwise, there exists a sequence $\varepsilon_n \downarrow 0$ such that
\begin{align*}
\lim_{n \rightarrow \infty}\sup_{x \in \R^3}\int_{B_{R}(x)}|u_{\varepsilon_n}^i|^2 \mathrm{d}x = 0
\end{align*}
for any $R>0$. By the Lions Lemma, we get $\int_{\R^3}|u_{\varepsilon_n}^i|^t\mathrm{d}x \rightarrow 0$ for any $t \in (2,6)$.
Then by $(f_1)$, we derive $\int_{\R^3}F(u_{\varepsilon_n}^i)\mathrm{d}x \rightarrow 0$ and $\int_{\R^3}f(u_{\varepsilon_n}^i)u_{\varepsilon_n}^i\mathrm{d}x \rightarrow 0$.
Since $\gamma_{\varepsilon_n}^i \rightarrow m_{V_0}$, by $I_{\varepsilon_n}(u_{\varepsilon_n}^i) = \gamma_{\varepsilon_n}^i$ and $I_{\varepsilon_n}'(u_{\varepsilon_n}^i) = 0$, we have
\begin{align*}
&m_{V_0} + o_n(1)= \frac{1}{2}\|u_{\varepsilon_n}^i\|_{\varepsilon_n}^2+\frac{b}{4}\left(\int_{\R^3}|\nabla u_{\varepsilon_n}^i|^2 \mathrm{d} x\right)^2-\frac{1}{6}
\int_{\mathbb{R}^{3}}|u_{\varepsilon_n}^i|^6\mathrm{d}x,\notag\\
& \|u_{\varepsilon_n}^i\|_{\varepsilon_n}^2+b\left(\int_{\R^3}|\nabla u_{\varepsilon_n}^i|^2 \mathrm{d} x\right)^2- \int_{\mathbb{R}^{3}}|u_{\varepsilon_n}^i|^6\mathrm{d}x = o_n(1).
\end{align*}
So
\begin{align*}
&m_{V_0} \ge \frac{a}{3}\int_{\R^3}|\nabla u_{\varepsilon_n}^i|^2 \mathrm{d} x+\frac{b}{12}\left(\int_{\R^3}|\nabla u_{\varepsilon_n}^i|^2 \mathrm{d} x\right)^2+ o_n(1),\notag\\
& a\int_{\R^3}|\nabla u_{\varepsilon_n}^i|^2 \mathrm{d} x +b\left(\int_{\R^3}|\nabla u_{\varepsilon_n}^i|^2 \mathrm{d} x\right)^2 \le \int_{\mathbb{R}^{3}}|u_{\varepsilon_n}^i|^6\mathrm{d}x + o_n(1).
\end{align*}
Similar to the argument of (\ref{3.9})-(\ref{3.10}), we can derive
$m_{V_0} \ge \hat{c}$, a contradiction with Lemma $2.2$.

\ep

\bl\lab{Lemma 4.2} $\varepsilon x_\varepsilon^i$ is bounded in $\R^3$, $i=1$, $2$, \ldots $k$.
\el

\bp  Fix $i=1$, $2$, \ldots $k$. Assume to the contrary that there exists a sequence $\varepsilon_n \downarrow 0$ such that $\varepsilon_n x_{\varepsilon_n}^i \rightarrow \infty$. By $I_{\varepsilon_n}(u_{\varepsilon_n}^i) = \gamma_{\varepsilon_n}^i \rightarrow m_{V_0}$ and $I_{\varepsilon_n}'(u_{\varepsilon_n}^i) = 0$, we have
$m_{V_0}+o_n(1)=I_{\varepsilon_n}(u_{\varepsilon_n}^i) -\frac{1}{4} \left(I_{\varepsilon_n}'(u_{\varepsilon_n}^i),u_{\varepsilon_n}^i\right) \ge \frac{1}{4}\|u_{\varepsilon_n}^i\|_{\varepsilon_n}^2$.
Then $\|u_{\varepsilon_n}^i\|_{\varepsilon_n}$ is bounded. Since the embedding $H_{\varepsilon_n} \hookrightarrow H$ is continuous, we know $\{u_{\varepsilon_n}^i\}$ is bounded in $H$. Set
$v_{\varepsilon_n}^i=u_{\varepsilon_n}^i(.+x_{\varepsilon_n}^i)$. Then by Lemma $4.1$, we get $\int_{B_{R_0}(0)}|v_{\varepsilon_n}^i|^2 \mathrm{d}x \ge \gamma_0$. So $v_{\varepsilon_n}^i \rightharpoonup v^i \ne 0$ weakly in $H$. Set
\begin{align*}
L_\varepsilon^i(u)=&\frac{1}{2} \int_{\R^3}\left(a|\nabla u|^2 +V(\varepsilon x+\varepsilon x_\varepsilon^i) |u|^2\right)\mathrm{d}x+\frac{b}{4}\left(\int_{\R^3}|\nabla u|^2 \mathrm{d} x\right)^2\\
&- \int_{\R^3}F(u) \mathrm{d}x-\frac{1}{6}\int_{\R^3}|u|^6\mathrm{d}x, \ \ u \in H_\varepsilon.
\end{align*}
By $I_{\varepsilon_n}(u_{\varepsilon_n}^i)  \rightarrow m_{V_0}$ and $I_{\varepsilon_n}'(u_{\varepsilon_n}^i) = 0$, we have
$L_{\varepsilon_n}(v_{\varepsilon_n}^i) \rightarrow m_{V_0}$ and $L_{\varepsilon_n}'(v_{\varepsilon_n}^i) = 0$.
By $\left(L_{\varepsilon_n}'(v_{\varepsilon_n}^i),v^i\right) = 0$, we get
\begin{align*}
&\int_{\R^3}\left(a \nabla v_{\varepsilon_n}^i \nabla v^i +V(\varepsilon_n x+\varepsilon_n x_{\varepsilon_n}^i) v_{\varepsilon_n}^i v^i \right)\mathrm{d}x+b\int_{\R^3}|\nabla v_{\varepsilon_n}^i|^2  \mathrm{d} x\int_{\R^3}\nabla v_{\varepsilon_n}^i \nabla v^i \mathrm{d} x\\
&= \int_{\R^3}f(v_{\varepsilon_n}^i) v^i\mathrm{d}x+\int_{\R^3}(v_{\varepsilon_n}^i)^5 v^i\mathrm{d}x.
\end{align*}
Since $v^i$ is non-negative, by $\varepsilon_n x_{\varepsilon_n}^i \rightarrow \infty$, $v_{\varepsilon_n}^i \rightharpoonup v^i$ weakly in $H$ and Fatou's Lemma, we derive
\begin{align}\label{4-1}
\|v^i\|_{V_\infty}^2+b\left(\int_{\R^3}|\nabla v^i|^2 \mathrm{d} x\right)^2 \le  \int_{\R^3}f(v^i)v^i \mathrm{d}x+ \int_{\R^3}|v^i|^6\mathrm{d}x.
\end{align}
By $(f_2)$, there exists a unique $t^i>0$ such that $t^i v^i \in M_\infty$. Similar to the argument of (\ref{2-1}), we get $t^i \le 1$.
By $L_{\varepsilon_n}(v_{\varepsilon_n}^i) \rightarrow m_{V_0}$ and $L_{\varepsilon_n}'(v_{\varepsilon_n}^i) \rightarrow 0$,
\begin{align}\label{4-3}
m_{V_0}=& L_{\varepsilon_n}(v_{\varepsilon_n}^i)-\frac{1}{4} \left(L_{\varepsilon_n}'(v_{\varepsilon_n}^i),v_{\varepsilon_n}^i\right) +o_n(1)\notag\\
= & \frac{1}{4} \int_{\R^3}\left(a|\nabla v_{\varepsilon_n}^i|^2 +V(\varepsilon_n x+\varepsilon_n x_{\varepsilon_n}^i) |v_{\varepsilon_n}^i|^2\right)\mathrm{d}x\notag\\
&+ \int_{\R^3}\left(\frac{1}{4} f(v_{\varepsilon_n}^i)v_{\varepsilon_n}^i - F(v_{\varepsilon_n}^i) \right)\mathrm{d}x+\frac{1}{12}\int_{\R^3}|v_{\varepsilon_n}^i|^6\mathrm{d}x+o_n(1).
\end{align}
Then by Fatou's Lemma and $t^i v^i \in M_\infty$ with $t^i \le 1$,
\begin{align*}
m_{V_0}  \ge & \frac{1}{4} \| v^i\|_{V_\infty}^2 + \int_{\R^3}\left(\frac{1}{4} f(v^i)v^i - F(v^i) \right)\mathrm{d}x+\frac{1}{12}\int_{\R^3}|v^i|^6\mathrm{d}x\notag\\
\ge & \frac{1}{4} \|t^i v^i\|_{V_\infty}^2 + \int_{\R^3}\left(\frac{1}{4} f(t^iv^i)t^iv^i - F(t^iv^i) \right)\mathrm{d}x+\frac{1}{12}\int_{\R^3}|t^iv^i|^6\mathrm{d}x\notag\\
=& I_{V_\infty}(t^i v^i)-\frac{1}{4} \left(I_{V_\infty}'(t^i v^i),t^iv^i\right)=I_{V_\infty}(t^i v^i) \ge m_{V_\infty} \ge m_{V_0},
\end{align*}
from which we derive $t^i=1$ and $V_\infty=V_0$. By (\ref{4-3}), we have
\begin{align}\label{4-4}
m_{V_0} \ge  \frac{1}{4} \| v_{\varepsilon_n}^i\|_{V_0}^2+ \int_{\R^3}\left(\frac{1}{4} f(v_{\varepsilon_n}^i)v_{\varepsilon_n}^i - F(v_{\varepsilon_n}^i) \right)\mathrm{d}x+\frac{1}{12}\int_{\R^3}|v_{\varepsilon_n}^i|^6\mathrm{d}x+o_n(1).
\end{align}
Then by Fatou's Lemma, $t^i=1$ and $V_\infty=V_0$,
\begin{align}\label{4-5}
m_{V_0}  \ge & \frac{1}{4} \| v^i\|_{V_0}^2 + \int_{\R^3}\left(\frac{1}{4} f(v^i)v^i - F(v^i) \right)\mathrm{d}x+\frac{1}{12}\int_{\R^3}|v^i|^6\mathrm{d}x\notag\\
=& I_{V_\infty}(t^i v^i)-\frac{1}{4} \left(I_{V_\infty}'(t^i v^i),t^iv^i\right)=I_{V_\infty}(t^i v^i) \ge m_{V_\infty} = m_{V_0}.
\end{align}
Combining (\ref{4-4})-(\ref{4-5}), we derive $v_{\varepsilon_n}^i \rightarrow v^i$ in $H$.
By the continuity of $\beta$, we have $\beta(v_{\varepsilon_n}^i) \rightarrow \beta(v^i)$.
On the other hand, by $u_{\varepsilon_n}^i \in N_{\varepsilon_n}^i$, we have $\beta(u_{\varepsilon_n}^i) \in C_{\frac{l}{\varepsilon_n}}^i$,
that is, $\frac{x^i-l}{\varepsilon_n}<\beta(u_{\varepsilon_n}^i)<\frac{x^i+l}{\varepsilon_n}$. Since $\beta(v_{\varepsilon_n}^i)=\beta(u_{\varepsilon_n}^i)-x_{\varepsilon_n}^i$, we get
\begin{align*}
\frac{x^i-l-\varepsilon_n x_{\varepsilon_n}^i}{\varepsilon_n}<\beta(v_{\varepsilon_n}^i)<\frac{x^i+l-\varepsilon_n x_{\varepsilon_n}^i}{\varepsilon_n},
\end{align*}
from which we derive $|\beta(v_{\varepsilon_n}^i)| \rightarrow \infty$, a contradiction with $\beta(v_{\varepsilon_n}^i) \rightarrow \beta(v^i)$.
So $\varepsilon x_\varepsilon^i$ is bounded in $\R^3$.

\ep

\bl\lab{Lemma 4.3} $\varepsilon x_\varepsilon^i \rightarrow x^i$ as $\varepsilon \rightarrow 0$, $i=1$, $2$, \ldots $k$.
\el

\bp  Fix $i=1$, $2$, \ldots $k$.
Since $\varepsilon x_\varepsilon^i$ is bounded in $\R^3$,
we can assume $\varepsilon x_\varepsilon^i \rightarrow x_0^i$ as $\varepsilon \rightarrow 0$.
Note that $I_{\varepsilon}(u_{\varepsilon}^i) = \gamma_{\varepsilon}^i \rightarrow m_{V_0}$ and $I_{\varepsilon}'(u_{\varepsilon}^i) = 0$.
Set $v_{\varepsilon}^i=u_{\varepsilon}^i(.+x_{\varepsilon}^i)$. By Lemma $4.1$, we have $v_{\varepsilon}^i \rightharpoonup v^i \ne 0$ weakly in $H$.
Since $\varepsilon x_\varepsilon^i \rightarrow x_0^i$ as $\varepsilon \rightarrow 0$, similar to the argument of Lemma $4.2$, we derive there exists $t^i \le 1$ such that $t^i v^i \in M_{V(x_0^i)}$.
By $I_{\varepsilon}(u_{\varepsilon}^i)  \rightarrow m_{V_0}$ and $I_{\varepsilon}'(u_{\varepsilon}^i) = 0$,
\begin{align}\label{4-6}
m_{V_0}=& I_{\varepsilon}(v_{\varepsilon}^i)-\frac{1}{4} \left(I_{\varepsilon}'(v_{\varepsilon}^i),v_{\varepsilon}^i\right) +o_\varepsilon(1)\notag\\
= & \frac{1}{4} \int_{\R^3}\left(a|\nabla v_{\varepsilon}^i|^2 +V(\varepsilon x+\varepsilon x_\varepsilon^i) |v_{\varepsilon}^i|^2\right)\mathrm{d}x+ \int_{\R^3}\left(\frac{1}{4} f(v_{\varepsilon}^i)v_{\varepsilon}^i - F(v_{\varepsilon}^i) \right)\mathrm{d}x\notag\\
&+\frac{1}{12}\int_{\R^3}|v_{\varepsilon}^i|^6\mathrm{d}x+o_\varepsilon(1).
\end{align}
Using Fatou's Lemma,
\begin{align}\label{4-7}
m_{V_0} & \ge \frac{1}{4} \int_{\R^3}\left(a|\nabla v^i|^2 +V(x_0^i) |v^i|^2\right)\mathrm{d}x+ \int_{\R^3}\left(\frac{1}{4} f(v^i)v^i - F(v^i) \right)\mathrm{d}x\notag\\
& \quad+\frac{1}{12}\int_{\R^3}|v^i|^6\mathrm{d}x\notag\\
&\ge I_{V(x_0^i)}(t^i v^i)-\frac{1}{4} \left(I_{V(x_0^i)}'(t^i v^i),t^i v^i\right)= I_{V(x_0^i)}(t^i v^i) \ge m_{V(x_0^i)} \ge m_{V_0},
\end{align}
in view of $t^i v^i \in M_{V(x_0^i)}$ with $t^i \le 1$. Then $t^i=1$ and $V(x_0^i)=V_0$. By (\ref{4-6}),
\begin{align*}
m_{V_0}\ge & \frac{1}{4} \|v_{\varepsilon}^i\|_{V_0}^2 + \int_{\R^3}\left(\frac{1}{4} f(v_{\varepsilon}^i)v_{\varepsilon}^i - F(v_{\varepsilon}^i) \right)\mathrm{d}x+\frac{1}{12}\int_{\R^3}|v_{\varepsilon}^i|^6\mathrm{d}x+o_\varepsilon(1).
\end{align*}
So by Fatou's Lemma,
\begin{align*}
m_{V_0} & \ge \frac{1}{4} \| v^i\|_{V_0}^2 + \int_{\R^3}\left(\frac{1}{4} f(v^i)v^i - F(v^i) \right)\mathrm{d}x+\frac{1}{12}\int_{\R^3}|v^i|^6\mathrm{d}x \\
&=I_{V(x_0^i)}(t^i v^i)-\frac{1}{4} \left(I_{V(x_0^i)}'(t^i v^i),t^i v^i\right)= I_{V(x_0^i)}(t^i v^i) \ge m_{V(x_0^i)} = m_{V_0},
\end{align*}
from which we derive $v_{\varepsilon}^i \rightarrow v^i$ in $H$. By the continuity of $\beta$, we get
$\beta(v_{\varepsilon}^i) \rightarrow \beta(v^i)$.
Together with $\beta(v_{\varepsilon}^i)=\beta(u_{\varepsilon}^i)-x_{\varepsilon}^i$ and $\beta(u_{\varepsilon}^i) \in C_{\frac{l}{\varepsilon}}^i$,
we get $\varepsilon x_\varepsilon^i \in C_{l}^i$ for $\varepsilon>0$ small. Then by $V(\varepsilon x_\varepsilon^i) \rightarrow V(x_0^i)=V_0$ as $\varepsilon \rightarrow 0$,
we derive $\varepsilon x_\varepsilon^i \rightarrow x^i$ as $\varepsilon \rightarrow 0$.

\ep

\bl\lab{Lemma 4.4} Fix $i=1$, $2$, \ldots $k$. Then $u_\varepsilon^i$ possesses a maximum $y_\varepsilon^i \in \R^3$ satisfying $V(\varepsilon y_\varepsilon^i) \rightarrow V(x^i)$ as $\varepsilon \rightarrow 0$. Moreover, there exist $C_0^i$, $c_0^i>0$ such that
\begin{align*}
u_\varepsilon^i(x) \le C_0^i \exp\left(-c_0^i |x-y_\varepsilon^i|\right)
\end{align*}
for $\varepsilon \in (0,\varepsilon_0)$ and $x \in \R^3$.
\el

\bp Let $v_{\varepsilon}^i=u_{\varepsilon}^i(.+x_{\varepsilon}^i)$.
Then $v_{\varepsilon}^i$ is the solution of
\begin{align}\label{4.6}
-\left( a+b \int_{\R^3}|\nabla v_{\varepsilon}^i|^2 \mathrm{d} x \right)\Delta v_\varepsilon^i +V(\varepsilon x+\varepsilon x_\varepsilon^i) v_\varepsilon^i
= f(v_{\varepsilon}^i) + (v_{\varepsilon}^i)^5 \ \ {\rm in } \ \  \R^3.
\end{align}
Moreover, by the argument of Lemma $4.3$, we know $v_{\varepsilon}^i \rightarrow v^i \ne 0$ in $H$ and $\varepsilon x_\varepsilon^i \rightarrow x^i$ as $\varepsilon \rightarrow 0$. Since $v_{\varepsilon}^i \rightarrow v^i$ in $H$, similar to the argument of Lemma $4.5$ in \cite{HZ}, we can derive $v_{\varepsilon}^i \in L^\infty(\R^3)$ and there exists $C>0$ independent of $\varepsilon \in (0,\varepsilon_0)$ such that $\|v_{\varepsilon}^i\|_\infty \le C$.
Moreover, $\lim_{|x| \rightarrow \infty}v_\varepsilon^i(x) =0$ uniformly for $\varepsilon \in (0,\varepsilon_0)$. Then by the elliptic estimate, there exist $\tilde{C}_0^i$, $c_0^i>0$ such that
\begin{align*}
v_\varepsilon^i(x) \le \tilde{C}_0^i \exp\left(-c_0^i |x|\right)
\end{align*}
uniformly for $\varepsilon \in (0,\varepsilon_0)$. Since the proof is standard, we omit it here.
Now we claim there exists $\varrho>0$ such that $\|v_{\varepsilon}^i\|_\infty \ge \varrho$ uniformly for $\varepsilon \in (0,\varepsilon_0)$. Otherwise, we have
$\|v_{\varepsilon}^i\|_\infty \rightarrow 0$ as $\varepsilon \rightarrow 0$. By (\ref{4.6}) and $V(\varepsilon x+\varepsilon x_\varepsilon^i) \ge V_0$, we get
\begin{align*}
\|v_{\varepsilon}^i\|_{V_0}^2 \le \int_{\R^3}\left(f(v_{\varepsilon}^i)v_{\varepsilon}^i + |v_{\varepsilon}^i|^6\right) \mathrm{d}x.
\end{align*}
From $(f_1)$, we have
$|f(v_{\varepsilon}^i)v_{\varepsilon}^i| + |v_{\varepsilon}^i|^6 \le \frac{V_0}{2}|v_{\varepsilon}^i|^2+C_{\frac{V_0}{2}}|v_{\varepsilon}^i|^6$,
where $C_{\frac{V_0}{2}}$ is a positive constant. Then
\begin{align*}
\|v_{\varepsilon}^i\|_{V_0}^2 \le 2 C_{\frac{V_0}{2}}\int_{\R^3}|v_{\varepsilon}^i|^6 \mathrm{d}x \le  2 C_{\frac{V_0}{2}}\|v_{\varepsilon}^i\|_\infty^4\int_{\R^3}|v_{\varepsilon}^i|^2 \mathrm{d}x \rightarrow 0
\end{align*}
as $\varepsilon \rightarrow 0$, a contradiction with $v_{\varepsilon}^i \rightarrow v^i \ne 0$ in $H$.

Let $z_\varepsilon^i \in \R^3$ be a maximum of $v_{\varepsilon}^i$, we have $\|v_{\varepsilon}^i(z_\varepsilon^i)\|_\infty \ge \varrho$.
Then by $\lim_{|x| \rightarrow \infty}v_{\varepsilon}^i(x)=0$ uniformly for $\varepsilon$, we derive there exists $M_0^i>0$ independent of $\varepsilon$ such that $|z_\varepsilon^i| \le M_0^i$. Since $z_\varepsilon^i$ is a maximum of $v_{\varepsilon}^i$, by $v_{\varepsilon}^i=u_{\varepsilon}^i(.+x_{\varepsilon}^i)$,
we know $x_\varepsilon^i+ z_\varepsilon^i$ is a maximum of $u_{\varepsilon}^i$. Let $y_\varepsilon^i=x_\varepsilon^i+ z_\varepsilon^i$.
Then by $\varepsilon x_\varepsilon^i \rightarrow x^i$ and $|z_\varepsilon^i| \le M_0^i$, we have $\varepsilon y_\varepsilon^i \rightarrow x^i$. Moreover, by $v_\varepsilon^i(x) \le \tilde{C}_0^i \exp\left(-c_0^i |x|\right)$ and $|z_\varepsilon^i| \le M_0^i$, we derive
\begin{align*}
u_{\varepsilon}^i(x) = v_{\varepsilon}^i(.-x_{\varepsilon}^i) \le \tilde{C}_0^i \exp\left(-c_0^i |x-x_\varepsilon^i|\right)=&\tilde{C}^i_0 \exp\left(-c^i_0 |x-y_\varepsilon^i+z_\varepsilon^i|\right)\\
\le& C_0^i \exp\left(-c_0^i |x-y_\varepsilon^i|\right).
\end{align*}

\ep

\vskip0.1in

\noindent\textit{\bf Proof of  Theorem  $1.1$ } By Lemmas $3.6$ and $4.4$, we know (\ref{2.1}) admits at least $k$ different positive solutions $u_\varepsilon^i$, $i=1$, $2$, \ldots $k$. Moreover, there exist $C_0^i$, $c_0^i>0$ such that
\begin{align*}
u_\varepsilon^i(x) \le C_0^i \exp\left(-c_0^i |x-y_\varepsilon^i|\right)
\end{align*}
for $\varepsilon \in (0,\varepsilon_0)$ and $x \in \R^3$. Let $v_\varepsilon^i=u_\varepsilon^i (\frac{.}{\varepsilon})$ and $z_\varepsilon^i=\varepsilon y_\varepsilon^i$.
Then $v_\varepsilon^i$ is the positive solution of (\ref{1.1})
and Theorem $1.1$ is proved.
 \hfill $\square$

\end{document}